%Musterdokument: doc.tex
\input amstex
\documentstyle{amsppt}
\magnification=1200
%\input avier %Anweisungen fuer dina4-Format
%Es folgt die "Praeambel", falls erforderlich
\overfullrule=0pt
\NoRunningHeads
%
%\input amacro1
%%%Makros-Dokument%%%allgemeine Makros
%%%%%%%%%%%%%%%%%%%%%%%%%%%%%%%%%%%%%%%%%%%%%
%Blackbold Buchstaben C,F,H,K,N,P,R,Q,Z: z.B. \C
%%%%%%%%%%%%%%%%%%%%%%%%%%%%%%%%%%%%%%%%%%%%%
\def\C{\Bbb C}
\def\F{\Bbb F}
\def\H{\Bbb H}
\def\K{\Bbb K}
\def\N{\Bbb N}

\def\P{\Bbb P}
\def\R{\Bbb R}

\def\Q{\Bbb Q}
\def\Z{\Bbb Z}
%%%%%%%%%%%%%%%%%%%%%%%%%%%%%%%%%%%%%%%%%%%%%%

%
\topmatter
\title
Hakenness and $b_1$
\endtitle
%
%Titel wird fortlaufende Kopfzeile auf den ungeraden Seiten, wenn %nicht
%anders abgegeben
%\rightheadtext{...}
\author
Alexander Reznikov 
\endauthor
%
%
%Autor wird Kopfzeile auf den geraden Seiten, wenn nicht anders angegeben
%\leftheadtext{...}
%\affil 
%Was bedeutet affill?
%\endaffil
%
%
%\address
%Adresse
%\endaddress
%
%
%\curraddr
%derzeitige Adresse
%\endcurraddr
%%
%\email
%E-Mail-Adresse
%\endemail
%
%
\date
September 1, 1997
\enddate

\date
Corrected version January, 1998
\enddate
%
%
%\dedicatory
%Widmung
%\enddedicatory
%
%
%\thanks
%Dank
%\endthanks
%
%
%\translator
%Uebersetzer
%\endtranslator
%
%
%\keywords
%Schluesselwoerter
%\endkeywords
%
%
%\subjclass
%Klassifizierung
%\endsubjclass
%
%
%\abstract
%Abstract
%\endabstract
%
%\input toc...
%Muster fuer Inhaltsverzeichnis: toc.tex
%\toc
%oder
%\toc\nofrills{Eigene Ueberschrift}
%\widestnumber\specialhead{...}
%\widestnumber\head{...}
%\widestnumber\subhead{...}
%\widestnumber\subsubhead{...}
%\head{A. Statements of the results} \endhead 
%
%\endtoc

\endtopmatter
\document
%\hyphenation{...}

\roster
\item"A." Statements of the main results
\item"B." Introduction

\item"   "
\item"1." Review of Thurston's hyperbolization theorem and quasi-Fuchsian representations 
\item"2." Review of Thurston's Surgery Theorem
\item"3." Review of Freedman's surfaces and Cooper-Long-Reid extension of Freedman's work
\item"4." Proof of the Surgery Theorem
\item"5." A fundamental trihotomy for
$\pi_1$-
injective surfaces
\item"6." Proof of the double coset theorem
\item"7." Quasi-Fuchsian knots and analytic proof of Freedman-Cooper-Long theorem
\item" 8."Freedman's surfaces give boundary groups, analysis $\grave a$ la Kerckhoff-Thurston
\item"9 ." Digression to Arithmetic Topology: Freedman's surface and generalized cyclotomic units
\item"10." Jaco's theorem and its corollaries
\item"11."
$p$-cycles 
\item"12." Homology domination
\item"13." Tight knots 
\item"14." Appendix: The language of arithmetic topology
\item"15." Appendix: Hyperbolization of Heegard splitted manifolds and bounded width in the mapping class group
\item"16." Appendix: Euler class and free generation
\item"17." Appendix: Homology of ramified coverings
\endroster

\head{A. Statements of the main results}
\endhead
We start with a new proof of the following unpublished theorem of Thurston:
\demo{A1. Theorem}
(Surgery Theorem). Let $M$ be a compact 3-manifold with boundary
$\partial M$
a union of $m$ tori. Assume $M$ is irreducible and atoroidal. Assume
that there is a finite covering $N$ of $M$, containing   an embedded
incompressible surface $S$ with no essential annuli joining $S$ and
$\partial N$.
Then for all parameters
$(m_i,n_i)$
of Dehn surgery with
$\underset i \to \min \sqrt{m^2_i+n^2_i} >> 1$,
the closed manifold
$\hat M(m_i,n_i)$
contains a
$\pi_1$
-injective surface.
\enddemo

\demo{A2. Double coset Theorem}
Let $M$ be an atoroidal virtually Haken manifold and let
$\pi_1(S)$
be a fundamental group of a virtually embedded surface. Then the double coset space
$\pi_1(S)\setminus \pi_1(M)/\pi_1(S)$
is infinite.
\enddemo

\demo{A3. Theorem}
(Freedman-Cooper-Long Theorem for quasi-Fuchsian knots). Let
$K \subset \Sigma$
be a quasi-Fuchsian knot. Let $S$ be an incompressible Seifert surface in
$\Sigma \setminus N(K)$.
Let
$m \in \pi_1(\Sigma \setminus K)$
be a meridian.
Then

$$m^n\pi_1(S)m^{-n} \cap \pi_1(S) =\{\text{cyclic group generated by [K]}\}$$
for
$n >> 1$.
\enddemo

\demo{A4. Theorem}
Let
$\Sigma \hookrightarrow M$
be a nonorientable surface. Assume
$H_1(\Sigma, \F_2) \rightarrow H_1(M,\F_2)$
is not surjective. Then $M$ has virtually positive
$b_1$.
\enddemo

\demo{A5. Theorem}
(Hakenness of loose ramified coverings) Let $M$ be a homology sphere with a link
$K = K_1 \cup \ldots \cup K_s$.
Let 
$N \rightarrow M$
be a $p$-covering, ramified along $L$. Let
$F_1$
be an oriented Seifert surface of
$K_1$,
transversal to
$K_i, i\ge 2$.
Let
$P = \sum^s_{i=2} \sharp (K_i \cap F_1)$.
Suppose for some
$i \ge 2$,
link
$(K_i,K_1) \ne 0(p)$.
If for some
$q \ne p$,

$$\dim H_1(N,\F_q) \ge 4-p(\chi(F)-2)+(p-1)P,$$
then $N$ is virtually Haken.
\enddemo

\demo{A6. Theorem}
(Ramified coverings of tight knots are Haken).Let
$K \subset \Sigma$
be a tight knot in a rational homology sphere. Let
$\tau : \Sigma_m \rightarrow \Sigma $
be a  ramified covering. Then $\Sigma_m$ is Haken or reducible.
\enddemo

\demo{A7. Corollary}Let
$K \subset \Sigma$
be any knot which is not fibered. Let $\hat \Sigma(m,n)$ be a
result of the Dehn surgery. Let $\hat \Sigma(m,n,p) $  be a
$p$-ramified covering of  $\hat \Sigma(m,n)$. Then for $m >> 1$ and
$p\geq 2$, $\hat \Sigma(m,n,p) $ is Haken or reducible.
\enddemo

%demo{A12. Theorem}
%Poincar$\acute e$ conjecture fails for number fields). There exists a number field
%K \ne \Q$
%uch that any unramified extension
%L \supset K$
%as degree 1.
%enddemo

%demo{A13. Conjecture}
%Bounded width in
%\text{Map}_g$).
%here are finitely many Dehn twists
%\lambda_i$
%n
%\text{Map}_g$
%uch that any double coset in

%$\text{Map}^0_g \setminus \text{Map}_g / \text{Map}^0_g$$
%as a representative of the form
%\prod^N_{i=1} \lambda^{m_i}_{j_i}$,
%or some fixed
%N \in \N$.
%nddemo

%demo{A14. Theorem}
%ssume A13. Then the Thurston's conjecture on hyperbolization of Heegrad splittings is true.
% $ddemo

\head{B. Introduction} \endhead

$   $\newline
%And you,no, really, Ivan, \newline
%You cannot argue that not all of your friends \newline
%Are impeccably dressed. \newline
%And some of them don't consume \newline
%Beverages of questionable quality. \newline
%In early ours of the morning.

This is a short summary of my work on some  open questions of
3-dimensional topology. 
%This part of mathematics has been terra icognita for most of mathematicians fo%r decades, developing very peculiar methods with essentially no interaction wi%th the rest of the current research. The situation has started to change when %Marden, J\o rgens
en, Riley, Gromov, and, principally, Thurston, brought a cycl%on of new ideas and methods into the subject. A very important discovery of Ba%ss, developed by Culler-Shalen, Morgan-Shalen and Culler-Gordon-Luecke-Shalen %brought in the light the role of re
presentation varieties and induced actions %on Tits buildings (trees in the case under consideration). These actions produ%ce incompressible surfaces (unfortunately, with boundary) in link complements,% which sometimes survive in the Dehn fillings [CGLS].

Since Haken and Waldhausen the most important question has been:
\roster
\item"(*)" How many closed 3-manifolds are Haken? \newline
The most optimistic conjecture is the following:
\item"(**)" Any irreducible manifold $M$ with infinite
$\pi_1$
is virtually Haken. \newline
The algebraic counterpart of (**) is
\item"(***)" Any $\pi_1(M)$ as above contains a surface group. \newline
And, finally, the strongest form of (**) is
\item"(****)" Any $M$ as above has virtually positive
$b_1$.
\endroster

Haken manifolds are "understood'' [Hem], so in case these questions are resolved positively, we would have a final conceptual picture of 3-manifolds. On the geometric side, this would resolve the hyperbolization conjecture positively for atoroidal 3-manif
olds with infinite
$\pi_1$
by an immediate application of Thurston's Hyperbolization Theorem for Haken 3-manifolds and a recent theorem of Gabai.

Three-manifolds can be constructed "practically'', 
%to the demand of research financing offices, and displayed in 5-space
by means of one of three competing constructions:

\demo{I. Surgery}
Take a three-manifold $M$ with boundary a union of tori
$T_1, \ldots T_m$
(a link complement in
$S^3$
is good), and glue in solid tori
$D^2 \times S^1$
by identifying
$\partial (D^2 \times S^1)$
with
$T_i$.
This depends on parameters
$(m_i,n_i)$,
once the basis for
$H_1(T_i,\Z)$
has been chosen.
\enddemo

\demo{II. Ramified covering}
Take a rational homology sphere
$\Sigma$
and a link
$L = \cup K_i$
in
$\Sigma$
and choose a homomorphism
$H_1(\Sigma \setminus L) \overset \psi \to \longrightarrow \Z$.
Form covering of
$\Sigma \setminus L$
corresponding to
$m \cdot ßZ$
in
$\Z$
and glue back the solid tori. This given a covering
$\Sigma_m \rightarrow \Sigma$
ramified along $L$. This depends on $m$ (and
$\psi$).
\enddemo

\demo{III. Heegard splitting}
Glue two handlebodies
$R^g_0, R^g_1$
by a homeomorphism
$\varphi : \partial R_0 \rightarrow \partial R_1$.
This depends on
$[\varphi] \in \text{Map}^0_g \setminus \text{Map}_g / \text{Map}^0_g$,
an element in the double coset space of the mapping class group
$\text{Map}_g$
factored by the actions of the subgroup of extendable homeomorphisms. 

There is a plenty of connections between these construction. Heegrad
splittings appear as surgeries, as explained in Appendix, section 22,
and ramified covering along knots are Heegard splittings unless some
compressions of Freedman's surfaces are nontrivial incompressible
surfaces , as explained in section 3.

As a part of this hyperbolization program Thurston proved the following three great theorems.
\enddemo

\demo{Theorem}
(Thurston) Haken implies hyperbolic. If $M$ is a closed Haken atoroidal manifold, then $M$ is hyperbolic.
\enddemo

\demo{Theorem}
(Thurston) Big surgeries are hyperbolic. If $M$ is an atoroidal manifold with boundary a union of tori, and
$\underset i \to \min \sqrt{m^2_i+n^2_i} > N$,
a fixed number then
$\hat M(m_i,n_i)$
is hyperbolic.
\enddemo

\demo{Theorem} 
(Thurston) Big ramified coverings are hyperbolic. If
$\Sigma$
is a rational homology sphere, $L$ a hyperbolic link in
$\Sigma, \psi : H_1(\Sigma \setminus L) \rightarrow \Z$
a homomorphism essential on meridians, then for
$m \ge N$,
a fixed number,
$\Sigma_m$
is hyperbolic.
\enddemo

Now, we have:

\demo{Theorem}
(Freedman-Cooper-Long) Big ramified coverings are Haken. If
$K \subset \Sigma$
is not a fiber knot than for $m$ big enough,
$\Sigma_m$
is Haken. That means Thurston's third theorem follows from the first.
\enddemo

For some surgeries, exactly those for which a big cyclic homology group is produced, this trivially implies that they are virtually Haken, see [CL]. An absolutely central question is:

$   $\newline
$(*^5)$ 
For $M$ an atoroidal manifold with a boundary a union of tori, do big surgeries contain a
$\pi_1$
-injective surface?

$   $ \newline
Our first result, theorem A1 is that the answer is {\bf yes}, if a
finite covering of $M$ contains an embedded incompressible closed
surface with no essential annuli joining it to the boundary. {\bf Any}
$M$ which does not fiber has a
finite covering $N$ which contains an  incompressible closed surface
by Freedman, but Freedman's surfaces are not good for us, as they
exactly {\bf do} have such annuli. Theorem A1 is not new: it appears
in [AR] without proof and is attributed there to Thurston. The
argument of Thurston, which we learned from J.H.Rubinstein, is
completely different from ours. Thurston uses the existence of metrics of non-positive curvature on Dehn surgeried manifolds, which are unchanged on the thick part of a cusped manifold, proved by Gromov-Thurston, where as we  work with quasi-Fuchsian defo
rmations in spirit of [Br].

Another central question, studied in many disguises by Jaco, is the following. Suppose we now already that $M$ contains a 
$\pi_1$-injective
surface. Is it embedded is some finite covering of $M$? In other words, is $M$ virtually Haken? A very fine theorem of Jaco [SW] states:

\demo{Theorem}
(Jaco) If
$\pi_1 (S) \hookrightarrow \pi_1(M)$
is an inclusion and if
$\pi_1(S)$
is contained in infinitely many subgroups of finite index in
$\pi_1(M)$,
then $S$ is embedded in a finite covering.
\enddemo

Now, we have the following fact, due to Thurston [Mo]

\demo{Proposition 5.1}
(Basic trihotomy) Let $M$ be a closed hyperbolic 3-manifold and
$\pi_1(S) \hookrightarrow \pi_1(M)$
an inclusion. Let 
$\Lambda$
be a limit set of
$\pi_1(S)$
acting on
$S^2$
by restriction of uniformization representation. Then
\roster
\item"(i)" if
$\Lambda \ne S^2$
or a Jordan curve, then $S$ is not embedded in any finite covering $N$ of $M$;
\item"(ii)" if
$\Lambda = S^2$
and $S$ is embedded in $N$, then $N$ fibers over a circle with fiber $S$;
 \item"(iii)" if
$\Lambda$
is a Jordan curve and $S$ is embedded in $N$ then $S$ is not a fiber of a fibration
$N \rightarrow S^1$
(but $N$ may still fiber).
\endroster
\enddemo

As a corollary, one shows easily 
\demo{Proposition 5.2}
(1st criterion of virtual fibering) Let $M$ be a closed hyperbolic three-manifold and
$\pi_1(S) \hookrightarrow \pi_1(M)$
an inclusion. Then there exists a finite covering $N$, which fibers over a circle with fiber $S$, if and only if
\roster
\item"1)" 
$\pi_1(S)$
is contained in infinitely many finite index subgroups of
$\pi_1(M)$
\item"2)" $\pi_1(S)$
contains a virtually normal subgroup of
$\pi_1(M)$.
\endroster
\enddemo
\demo{ Proposition 6.4}
(2nd criterion of virtual fibering). Let $M$ be a hyperbolic 3-manifold,
$\pi_1(S) \hookrightarrow \pi_1(M)$
an injective surface. Then $S$ is a fiber of a fibration for some finite covering $N$ of $M$ iff
\roster
\item"(i)" $\pi_1(S)$
is in infintely many subgroups of finite index and
\item"(ii)" for any
$x \in \pi_1(M)$
there exists
$N(x)$
such that for
$m \ge 1$

$$x^{mN} \pi_1(S) x^{-mN} \cap \pi_1(S) \ne \{1\}$$
\endroster
\enddemo
We will deduce a following result, which gives an answer to a question, posed to us by Michel Boileau:
\demo{A2. Double coset Theorem}
Let $M$ be an atoroidal virtually Haken manifold and let
$\pi_1(S)$
be a fundamental group of a virtually embedded surface. Then the double coset space
$\pi_1(S)\setminus \pi_1(M)/\pi_1(S)$
is infinite.
\enddemo

\demo{Quasi-Fuchsian knots}
Let $M$ be an atoroidal manifold with boundary a torus $T$ and let
$K \subset T$
be a simple loop, which is homologically trivial in $M$. We say $K$ is a quasi-Fuchsian knot if there exists an incompressible Seifert surface $S$ for $K$ such that the restriction of the uniformization representation of $M$ on
$\pi_1(S)$
is quasi-Fuchsian.
\enddemo

\demo{\bf Fact}
Any knot which is not fibered, and not 2-virtually fibered, is quasi-Fuchsian.
\enddemo

Using this fact  we give a simple geometric proof of the main result of Freedman [FF], and Cooper-Long [CL], which says that two lifts of $S$ in cyclic coverings
$M_n, n>>1$,
of $M$ glue together to form a closed incompressible surface in $M_n$.
It is elementary to show [CL] that these surface stay incompressible in all non-longitudal Dehn surgeries of 
$M_n$.

We now turn to a more algebraic part of the paper.

In two preprints distributed a year ago [Re 2] [Re 3] I introduced a notion of $p$-cycles domination:

Let
$Q \hookrightarrow M$
be an embedded $p$-cycle (for
$p=2$,
this is just a non-orientable surface). Then we say that $Q$ dominates $M$ if the image of
$\pi_1(Q) \rightarrow \pi_1(M)$
is of finite index in $M$.

{\bf Any 3-manifold with rich fundamental group [Re 1] has a finite covering containing} {\bf a $p$-cycle, for any prime $p$.} See [Re 1]. {\bf If $\Sigma$
is a homology sphere, $K$ a knot in
$\Sigma$
and}
$p,q$
{\bf coprime numbers, than}
$\hat \Sigma (p,q)$
{\bf carries a} $p$ {\bf -cycle with genus equal to the genus of the
knot, and therefore independent on} $p$. The Domination theorem 11.6  states a non-dominating $p$-cycle carries an 
$\pi_1$-injective
surface in $M$. The Cooper-Long theorem that for $K$ not fibered
$\hat \Sigma (p,q)$
are virtually Haken is equivalent to the algebraic fact that the $p$-cycle mentioned above is not dominating for $p$ big enough.

For a group $G$, let 
$r(G)$
be the least number of generators for $G$.

\demo{Problem}
If $K$ is not a fibered knot in
$\Sigma$,
then
$r(\pi_1(\Sigma_m)) \underset{m \rightarrow \infty} \to \longrightarrow \infty$.
I explain in $\qquad$ why this implies Freedman's theorem that
$\pi_1(\Sigma \setminus K)$
contains surface groups and
$\Sigma \setminus K$
virtually contains an incompressible closed  surface.

A $p$-cycle is $q$-homology dominating, if
$H_1(Q,\F_q) \rightarrow H_1(M,\F_q)$
is surjective. For
$p = q=2$
we have the theorem A4, which says that if a nonorientable surface is not 2-homology dominating, then $M$ has a double unramified covering with positive
$b_1$.
In other words,
$\pi_1(M)$
maps onto
$C_2 * C_2$.\footnote{By a theorem of Rubinstein, $M$ itself is Haken}

A large variety of 3-manifolds with not-dominating $p$-cycles come from ramified coverings over "loose" links. By the latter I mean a link
$L = K_1 \cup \ldots \cup K_s$
in a homology sphere
$\Sigma$
such that
$K_j, j \ge 2$
are "much knotted" outside Seifert surface $S$ of
$K_1$. We fix $S$ and demand that the number of intersections of
$K_j$
with $S$ stays bounded, whereas 
$K_j$
are complicated outside $S$. If
$\hat \Sigma_p$
is a ramified covering over $L$ then a simple condition ensure that the preimage of $S$ in
$\hat \Sigma_p$
is an embedded non-dominating cycle
$\text{mod} \, p$, and
$\hat \Sigma_p$
is virtually Haken (the last statement uses a criterion of Shalen-Wagreich [SW]).
\enddemo

\demo{Tight knots}
A knot
$K \subset \Sigma$
is called tight, if there is a Seifert surface $S$ with
$\partial S =K$,
such that
$\pi_1(S) \rightarrow \pi_1(\Sigma)$
is injective (for a Seifert surface of minimal genus one always has
$\pi_1(S) \rightarrow \pi_1(\Sigma \setminus K)$
injective [BZ]). There are many tight knots. Take any free finitely generated group in
$\pi_1(M)$
(there are plenty of such for manifolds with rich fundamental group, [SW]) and realize it by an incompressible handlebody 
$D^g$
in $M$. Take an essential simple loop $\ell$ in
$\partial D$
which is in the kernel of
$H_1(\partial D) \approx \Z^{2g} \rightarrow \Z^g = H_1(D)$
and take an incompressible surface $F$ in $D$ which span this loop. One sees immediately that
$\ell$
is tight in $M$. {\bf Any knot which is not fibered, is tight in
sufficiently big ramified coverings}. This is an immediate consequence
of Freedman-Cooper-Long theorem.{\bf Any knot which is not fibered, is
tight in a $(m,n)$-surgery, if $m$ is big enough}.This is a resent
theorem of Boyer, Culler, Shalen and Zhang [BCSZ].  Any genus one knot in a hyperbolic manifold $M$, whose fundamental group is not virtually generated by 2 elements, is tight or there are
$\pi_1$-injective
closed surfaces in $M$.

Let
$\Sigma_m \rightarrow \Sigma$
be a hyperbolic ramified covering over a tight knot. Our Theorem A6
asserts that $\Sigma_m $ is Haken.
We do not assume $m$ is big
$(m=2$
is fine). 
The proof of the theorem A6 uses a technique of {\bf collapse}. The first appeared in [Re ], which formed a part of this work.

\demo{Arithmetic topology}
As I already wrote in [Re 1] the underlying philosophy of my approach to 3-manifolds is the idea that 3-manifold topology is parallel to algebraic number theory. The manuscript [KR] based on a lecture of Kapranov and myself in MPI in 1996 contained the di
ctionary joining the two branches of mathematics together. With a kind approval of Kapranov, I am now publishing this dictionary, in Appendix 14.
\enddemo

\demo{Bounded width and Thurston's conjecture on hyperbolization of big Heegard splittings}
An exciting recent discovery of bounded width asserts that for arithmetic lattices
$\Gamma$
in semisimple algebraic groups of rank
$\ge 2$,
there is a set
$x_1, \ldots x_n$
of generators such that any
$x \in \Gamma$
is represented as
$x^{n_1}_{i_1} \ldots x^{i_N}_{i_N}$
with $N$ a fixed number. In appendix 15  I explain how the Thurston conjecture follows from an algebraic conjecture of bounded width in the modular group
$\text{Map}_g$.
\enddemo

\demo{Euler class and free generation}
The appendix 16 contains a Fuchsian model for doubling procedure which
leads from an incompressible Seifert surface to Freedman surface.The
main result of this appendix , Theorem 16.0, gives sufficient conditions for
$n$ elements of $PSL_2(\R)$ to generate a discrete group.
\enddemo

\demo{Homology of ramified coverings}
The appendix 17 contains a theorem A8, which describes the Galois module
$H_1(N)_{(p)}$
for a
$C_p$
-ramified covering $N$ over a link
$L \subset M$
with $s$ components.This theorem parallels the approach of [Re 1],
the Rough Classification Theorem.
\proclaim{Theorem A8}
Let
$N \rightarrow M$
be a covering of rational homology spheres ramified along a link with $s$ components. Then as a 
$C_p$-module,
$H_1(N)_{(p)}$
is a direct sum of exactly
$(s-1)$
open modules and some number of closed modules.
\endproclaim

 We prove also a three-manifold
version of a classical results of Gauss and Redei about the 2-part of
the ideal class group of quadratic fields, theorem A9.
\proclaim{Theorem A9}
Let
$p=2$
and
$s=2$,
so that the link is
$K_1 \cup K_2$
\roster
\item"(a)" if the linking number
$\text{link}\, (K_1,K_2)$
is odd, then
$H_1(M)_{(2)} = \Z/2\Z$.
\item"(b)" if the linking number
$\text{link}\, (K_1,K_2)$
is even, then either
$b_1(N) > 0$
or
$H_1(N)_{(2)} = \Z/2^{\ell}\Z$
with
$\ell \ge 2$.
\endroster
\endproclaim

\bf{Acknowledgements}
I would like to thank Oertel, Boileau, Boyer, Barge, Lott, Rubinstein, Potyagailo for useful discussions, MPI and Manuela Sarlette for expert typing. As well as [Re 1] this paper would never have been written without the support of Lucy Katz (Reznikov).

\bf{1. Review of Thurston's hyperbolization theorem and quasi-Fuchsian representations [Mo] [Th 2] [Mc 1] [Mc 2]}

\bf{1.1 Quasi-Fuchsian representations}

A discrete group
$\Gamma$
in
$PSL_2(\C)$
is called Quasi-Fuchsian if the limit set
$\Delta$
is a Jordan curve, and Fuchsian if
$\Delta$
is a circle.

Let
$\Gamma$
be a Fuchsian group leaving the upper half-plane  $U$ and the lower half-plane $L$ invariant. A Beltrami coefficient 
$\mu$
of norm $< 1$
supported in $U$ is invariant, if a natural action of
$\Gamma$
in
$\bar K \cdot K^{-1}$
leaves it invariant, where $K$ is the canonical bundle. Such
$\mu$
determines a quasi-conformal homeomorphism (normalized)
$f : \bar \C \rightarrow \bar \C$,
conjugating
$\Gamma$
to a quasi-Fuchsian representation. On the other hand, extend
$\mu$
to the lower half-plane by conjugation; let
$\nu$
be the resulted Beltrami coefficient. The correspondent quasiconformal map
$\psi$
will conjugate 
$\Gamma$
to another Fuchsian representation. We have therefore a correspondence:

$$\text{pair of points in the Teichm\"uller space of} \quad \Gamma \leftrightarrow \quad \text{quasiconformal representations}$$ 
intensively studied by Bers [B]. The map
$f \mid L$
above is schlicht holomorphic. The correspondence

$$\text{quasiconformal representation} \quad \rightarrow \quad \text{Schwartzian of} \quad f$$
is called the Bers embedding. Its image contains a ball in
$H^0(K^2)$
with
$L^1$
-norm and is contained in a ball [Gar]. The boundary of this image correspond to some representation
$\Gamma \rightarrow PSL_2(\C)$,
called boundary groups. If
$\Gamma = \pi_1(S)$,
a fundamental group of a closed surface, and
$\rho : \Gamma \rightarrow PSL_2(\C)$
is quasi-Fuchsian, then a neighbourhood of the convex core $M$ of
$\Cal H^3/\rho(\Gamma)$
is homeomorphic to
$S \times [0,1]$.
We see therefore that quasi-Fuchsian representations are parametrized by
$T(\partial M)$.
\enddemo

\demo{1.2 Deformation of Kleinian manifolds}
Let
$\Gamma$
be a finitely generated Kleinian group such that all components of
$S^2\setminus \Delta (\Gamma)$
are simply connected, and the quotient of the  neighbourhood of the convex hull of
$\Delta(\Gamma)$ by  $\Gamma$
has finite volume. The manifold
$\Cal H^3/\Gamma \cup (S^2\setminus \Delta (\Gamma))/\Gamma$
is called geometrically finite. We are interested to deform
$\Gamma$,
by quasiconformal homeomorphism of
$S^2$,
keeping these properties. A {\bf deformation theorem} [BK] asserts that, again, the space of such deformations is identified with
$T(\partial M)$.
Topologically, the thick part
$\{\text{injectivity radius} \quad (x) \ge \varepsilon\}$
of $M$ is compact with boundary consisting of
$\partial M$
and a number of tori, corresponding to cusps of $M$. We will not make any difference between
$T(\partial M)$
and the space of quasiconformal deformation of $M$.
\enddemo

\demo{1.3 The skinning map}
For any component $S$ of
$\partial M$,
the restriction of the representation to
$PSL_2(\C)$
on
$\pi_1(S)$
is quasi-Fuchsian. The representation of
$\pi_1(S)$
is another component of
$S^2\setminus \Delta (\pi_1(S))$
gives, by Bers a point in
$T(S)$.
One gets a map

$$T \undersetbrace \theta \to {(\partial M) \rightarrow \text{Kleinian groups for} \quad \pi_1(M) \rightarrow T} (\partial M)$$
called skinning map. It is holomorphic, therefore a contraction for Teichm\"uller metric by Royden's theorem.
\enddemo

\demo{1.4 The Thurston's main idea}
Let
$M = M_1 \cup M_2$
identified by homeomorphisms of some bounadry components. Suppose we realized
$M_i$
as geometrically finite Kleinian manifolds. We want to deform
$M_i$
such that the restriction of the representations on corresponding boundary components will correspond by an induced map in the Teichm\"uller space (it means induced by the gluing map). Then by Van Kampen we will get a representation for $M$. The skinning 
map followed by this latter isometry gives a holomorphic self-map in
$T(\partial M)$
and we want to show that this map has a fixed point {\bf inside}
$T(\partial M)$.
One first extends
$\theta$
to the closure of
$T(\partial M)$
in the representation variety

$$V = \text{Hom} (\pi_1(M), PSL_2 (\C))/PSL_2(\C)$$
and then shows that, if $M$ is atoroidal, the image of this extension is precompact in
$T(\partial M)$,
so
$\theta$
is uniformly conracting. See [Mo] for the outline of this proof.
\enddemo

\demo{1.5 McMullen's proof} This proof is exposed in details in the  recent work of
Otal [O].
\enddemo

\demo{1.6 Resum$\acute e$}
What we really will need from the above, are the following theorems, are immediate corollaries of the above discussion, cd. [Mo], [Mc 1] [Mc 2].
\enddemo

\demo{Theorem 1.6.1}
Let $M$ be an atoroidal manifold with boundary a union of tori. Let
$S \hookrightarrow M$
be an embedded incompressible surface. Suppose there are no essential annuli joining $S$ to the boundary. Then the restriction of the uniformization representation of $M$ on
$\pi_1(S)$
is quasi-Fuchsian.
\enddemo

Similarly:

\demo{Theorem 1.6.2}
Let $M$ be an atoroidal Haken manifold. Let
$S \hookrightarrow M$
be an incompressible surface, not a fiber of a fibration of $M$ over
$S^1$
. Then the restriction of the uniformization representation of $M$ on
$\pi_1(S)$
is quasi-Fuchsian.
\enddemo

\demo{Theorem 1.6.3}
Let $M$ be a  manifold of Betti number one with boundary a torus,
whose interior carries a hyperbolic metric of finite volume.Let $S$ be a Seifert surface for $M$ which is not a fiber of a fibration over $S^1$ and not a 2-virtual fiber.Then the restriction of the uniformization representation of $M$ on $S$ is quasi-Fuch
sian.
\enddemo

We also notice:

\demo{Theorem 1.6.4}
The space of quasi-Fuchsian representation of
$\pi_1(S)$,
a fundamental group of a closed surface, is open in the representation variety

$$\text{Hom} \, \pi_1(S), PSL_2(\C))/PSL_2(\C)$$
\enddemo

\bf{2. Review of Thurston's surgery theorem} 

\demo{2.1}
\enddemo
Let
$\Gamma$
be a countable group. We can endow the quasiprojective variety

$$\text{Hom}_{irr}(\Gamma, PSL_2(\C))/PSL_2(\C)$$
with the classical topology. If
$\rho_i$
is a family of irreducible representations, we say
$\rho_i$
converges algebraically to a representation
$\rho$
if the classes of
$\rho_i$
converge to the class of
$\rho$
in the classical topology. An important fact, due to Chaubaty and J\o rgensen, is that if
$Im \, \rho_i$
are discrete, so is
$\rho$.

Now let $M$ be an atoroidal manifold with boundary a union of tori
$T_i$.
Fix a basis
$(e_i,f_i)$
for
$H_1(T_i,\Z) \approx \Z \oplus \Z$.
For
$(m_i,n_i)$
coprime we denote
$\hat M (m_i, n_i)$
a result of the Dehn surgery which kills the element
$m_ie_i + n_if_i \in \pi_1(T_i)$.
\enddemo

\demo{Theorem 2.1}
\roster
\item"(a)" For
$\underset i \to \min \sqrt{m^2_i + n^2_i} >> 1, \hat M(m_i,n_i)$
is hyperbolic.
\item"(b)"
The composition of the surjective map
$\pi_1(M) \rightarrow \pi_1 (\hat M (m_i, n_i))$
and the uniformization representation of
$\hat M$
coverges algebraically to the uniformization representation of $M$ as
$\underset i \to \min \sqrt{m^2_i+n^2_i} \rightarrow \infty$
\endroster
The proof is  "elementary'' in the sense that it does not use PDE, like
the existence of quasiconformal maps with given Beltrami coefficient. It is somewhat relevant to the elementary proof of the last mentioned result given in [LV]. Hatcher and Thurston proved that for 2-bridge knots the surgeries
$\hat M$
are not Haken, though they are hyperbolic [HT]. On the other hand by
Theorem  2.03  of [CGLS] if $M$ contains an incompressible closed surface, then big surgeries
$\hat M$
are Haken. Our theorem A1 states that if a {\bf covering} of $M$
contains an incompressible surface with no essential annuli joining it to the boundary, then big surgeries
$\hat M$
contain a
$\pi_1$
-injective surface. There is little doubt that they are virtually Haken.
\enddemo

\head{3. Review of Freedman's surfaces and Cooper-Long extension of Freedman's work} \endhead

\demo{3.1}
Let $G$ be a group. Suppose
$x_i,y_i\mid^g_{i=1}$
are elements of $G$,
$z = \prod^g_{i=1} [x_i,y_i]$
and
$\sigma : G \rightarrow G$
an automorphism, fixing $z$.
\enddemo

\demo{Definition}
A Freedman surface
$Fr(x_i,y_i,\sigma)$
in
$K(G,1)$
is a map of
$\hat S^{2g}$( a closed surface of genus 2g)
to
$K(g,1)$
defined as follows: consider the amalgam

$$\pi_1(S^{2g}) = F_g \underset{\Z} \to * F_g$$
and map it to $G$ by sending the generators of the first copy
$F_g$
to
$x_i,y_i$
and generators of the second copy to
$\sigma(x_i), \sigma(y_i)$.
\enddemo

\demo{Important  example}
Let
$G = \pi_1(M)$,
a fundamental group of an irreducible rational homology sphere. Assume $G$ infinite, so $M$ is acyclic. Suppose for some prime
$p, H_1(M,\Z)_{(p)} \ne 0$.
Let
$w \in \pi_1(M)$
be an element with projects nontrivially to
$H_1(M,\Z)_{(p)}$.
Let
$p^n$
be an order of $w$, so that

$$z = w^{p^n} \in [\pi_1(M), \pi_1(M)].$$
Let
$x_i, y_i \in \pi_1(M)$
be such that
$\prod^n_{i=1} [x_i,y_i] = w^{p^n}$.
Then one forms a Freedman surface
$Fr(x_i, y_i, Ad \, w)$.
We will see in section $\qquad$ that sometimes this surface is 
$\pi_1$-injective in $M$.
\enddemo

\demo{Freedman's surface by Freedman-Freedman [FF]}
The following special case has been considered in [FF]. Let $M$ be a
3-manifold  having the first Betti number one with a boundary a torus $T$, and let $K$ be a longitude, i.e. an essential simple loop in $T$, homologically trivial in $M$. Let
$S^g$
be a Seifert surface for $K$. We can assume $S$ is
$\pi_1$-injective
by compression [BZ]. Let $m$ be another element of
$\pi_1(T)$.
Let
$x_i, y_i$
be canonical generators for
$\pi_1(S)$
so that
$\prod^g_{i=1} [x_i, y_i] = [K]$.
Then one looks at
$Fr (x_i, y_i, Ad \, m^{\kappa})$.
Obviously this surface lifts to an infinite cyclic covering
$M_{\infty}$
of $M$.
\enddemo

\demo{Theorem 3.1} 
Suppose $S$ is not a fiber of a fibration over
$S^1$.
\roster
\item"(a)" for
$\kappa >> 1$,
the image of
$\pi_1(Fr (x_i, y_i, Ad \, m^{\kappa}))$
in
$\pi_1(M_{\infty})$
is not free.
\item"(b)" in some finite cyclic covering of $M$ there is a closed embedded incompressible surface obtained by compression of
$Fr(x_i, y_i, Ad \, m^{\kappa})$
\endroster
\enddemo

\demo{Cooper-Long extension}
\enddemo

\demo{Theorem 3.2}
In condition of Theorem 3.1
\roster
\item"(a)" $Fr(x_i, y_i, Ad \, m^{\kappa})$
is
$\pi_1$
-injective for
$\kappa \ge 6g-2$
\item"(b)" all finite coverings 
$M_{\kappa}$
with
$\kappa \ge 6g-2$
contain a closed incompressible surface with an annulus joining it to the boundary torus
\item"(c)" this incompressible surface survives in all Dehn surgeries of
$M_{\kappa}$
except one, killing the longitude
\item"(d)" for
$(p,q)$
coprime and
$p \ge 6g-2$,
the Dehn surgery
$\hat M(p,q)$,
having a nontrivial $p$-part of
$H_1(\hat M,\Z)$,
has a $p$-cyclic covering which is Haken, or reducible.
\item"(e)" the ramified coverings
$N_{\kappa}$
of any Dehn surgery $N$ of $M$, except longitudal are Haken for
$k >> 1$.
\endroster
We notice that (b)-(e) are immediate consequences of (a). The proof in [FF] and [CL] uses hyerarchies; an approach of [CLR] is much more elementary, but it does not give an information, applicable to Dehn surgeries. We will be proving our version, Theorem
 A3, by analysis of the limit set in section
\enddemo

\head{4. Proof of the surgery theorem, big surgeries contain a
$\pi_1$-injective surface} \endhead

\demo{Proof of the Theorems 1.6.1, 1.6.2}
Cut $M$ along $S$ and let
$M_0$
be the resulted manifold. Then by Thurston,  there is a geometrically finite realization of
$M_0$
with incompressible quasi-Fuchsian boundary, such that the gluing problem has a solution. This exactly means that the restriction of the uniformization representation on
$\pi_1(S)$
is quasi-Fuchsian.
\enddemo

\demo{Proof of the Surgery Theorem}
By the virtue of Theorem 1.6.1, applied to $N$ the restriction of the uniformization representation of $M$ on
$\pi_1(S)$
is quasi-Fuchsian. By Thurston's Surgery Theorem 2.1 (a),
$\hat M(m_i,n_i)$
are hyperbolic for
$\underset i \to \min \sqrt{m^2_i+n^2_i} >> 1$.
Let
$\rho_i$
be the composition of the surjective map
$\pi_1(M) \rightarrow \pi_1(\hat M)$
and the uniformization representation of
$\hat M$.
As
$\underset i \to \min \sqrt{m^2_i + n^2_i} \rightarrow \infty$,
$\rho_i\mid S$
converge to a quasi-Fuchsian representation, by the above and Theorem 2.1 (b). Now, by Theorem 1.6.4 the set of quasi-Fuchsian representation is open in the representation variety of
$\pi_1(S)$.
Therefore for
$\underset i \to \min \sqrt{m^2_i+n^2_i}$
big enough,
$\rho_i$
is quasi-Fuchsian, therefore injective. Hence
$\pi_1(S) \rightarrow \pi_1(\hat M)$
is injective.
\enddemo

\head{5. A fundamental trihotomy for $\pi_1$-injective surfaces} \endhead

Let $M$ be a closed 3-manifold and
$S \rightarrow M$
a
$\pi_1$
-injective immersed surface. A classical question is: does there exists a finite covering $N$ of $M$ with
$\pi_1(N) \supset \pi_1(S)$
such that $S$ lifts an embedded surface in $M$? A theorem of Jaco [SW] answers this question affirmatively if
$\pi_1(S)$
lies in infinitely many subgroups of finite index. Now assume $M$ to be hyperbolic. One has

\demo{Proposition 5.1}
(Basic trihotomy) Let $M$ be a closed hyperbolic 3-manifold and
$\pi_1(S) \hookrightarrow \pi_1(M)$
an inclusion. Let 
$\Lambda$
be a limit set of
$\pi_1(S)$
acting on
$S^2$
by restriction of uniformization representation. Then
\roster
\item"(i)" if
$\Lambda \ne S^2$
or a Jordan curve, then $S$ is not embedded in any finite covering $N$ of $M$
\item"(ii)" if
$\Lambda = S^2$
and $S$ is embedded in $N$, then $N$ fibers over a circle with fiber $S$
 \item"(iii)" if
$\Lambda$
is a Jordan curve and $S$ is embedded in $N$ then $S$ is not a fiber of a fibration
$N \rightarrow S^1$
(but $N$ may still fiber)
\endroster
\enddemo

\demo{Proof}
Let $N$ be a finite covering of $M$ where $S$ is embedded. If $S$ is a fiber of a fibration over
$S^1$,
then
$\pi_1(S)$
is normal in
$\pi_1(M)$,
therefore
$\Delta (\pi_1(S))$
is
$\pi_1(M)$
-invariant, hence
$\Delta (\pi_1(S)) = S^2$.
If $S$ is not a fiber, then
$\pi_1(S)$
acts on
$S^2$
as a quasi-Fuchsian group, by Theorem 1.6.2, so
$\Delta (\pi_1(S))$
is a Jordan curve.
\enddemo

\demo{Proposition 5.2}
(1st criterion of virtual fibering) Let $M$ be a closed hyperbolic three-manifold and
$\pi_1(S) \hookrightarrow \pi_1(M)$
an inclusion. Then there exists a finite covering $N$, which fibers over a circle with fiber $S$, if and only if
\roster
\item"1)" 
$\pi_1(S)$
is contained in infinitely many finite index subgroups of
$\pi_1(M)$
\item"2)" $\pi_1(S)$
contains a virtually normal subgroup of
$\pi_1(M)$.
\endroster
\enddemo

\demo{Proof}
The "only if" part is obvious. Now, (i) implies by Jaco that $S$ is embedded is some finite covering $N$. On the other hand, if 
$N^{\prime}$
is another finite covering such that
$\pi_1(S)$
contains a nontrivial subgroup, normal in
$\pi_1(N)$,
then the limit set of this normal subgroup is all of
$S^2$
by above, so
$\Delta (\pi_1(S))$
is all of
$S^2$
and $S$ is a fiber of a fibration.
\enddemo

\head{6. Proof of the  double coset Theorem}\endhead
\demo{Proposition 6.1(Malnormality) }
Let $M$ be atoroidal Haken manifold,
$S \hookrightarrow M$
an embedded incompressible surface, not a fiber of a fibration over
$S^1$.
Then for a big set of elements (see below)
$x \in \pi_1(M)$
there exists
$N(x)$
such that for
$n \ge N(x)$

$$x^{-n}\pi_1(S)x^n \cap \pi_1(S) =1$$

\enddemo

\demo{Proof}By the Thurston's hyperbolization theorem we may, and will, consider $M$ as hyperbolic manifold. Let
$\Delta$
be
$\Delta (\pi_1(S))$,
a Jordan curve by 7.1. Define
$B \subset \pi_1(M)$
by
$x \in B$
if and only if neither of the two fixed points of $x$ is in
$\Delta$.

\demo{Lemma 6.2}
$B$ is not empty.
\enddemo

\demo{Proof}
Let
$C \subset S^2 \times S^2$
be a closure of the set of oriented pairs (attracting fixed point $(y)$), repelling fixed point $(y)$) taken over all
$y \in \pi_1(M)$. 
Then $C$ is closed and invariant, therefore $C$ is all of
$S^2 \times S^2$
by Sullivan's theorem [Su]. So $B$ is nonempty. Let
$y \in \pi_1(M)$
and let
$\Delta^{\prime} = y(\Delta)$.
Define
$B^{\prime}$
as a set of
$x \in \pi_1(M)$
such that neither of the fixed points of $x$ in
$\Delta \cup \Delta^{\prime}$.
\enddemo

\demo{Lemma 6.3}
For
$x \in B^{\prime}, x^ny^{-1} \in B$
for
$n >> 1$.
\enddemo

\demo{Proof}
The attracting fixed point of
$x^ny$
converges to that of $x$. The repelling fixed point of
$x^ny$
converges to the image under
$y^{-1}$
of the repelling fixed point of $x$.

For
$\Delta$
a closed set of
$S^2$
of measure $0$ we will call the set
$B(\Delta)$
defined as above {\bf a big set} of
$\pi_1(M)$.

Now let
$x \in B(\Delta (\pi_1(S))$.
Let $n$ be a positive integer and let
$y \in x^n\pi_1(S)x^{-n}$.
Then the fixed points of $y$ lie in
$x^{-n} (\Delta (\pi_1(S))$.
So if
$y \in x^n \pi_1(S)x^{-n} \cap \pi_1(S)$
then the fixed points of $y$ lie in
$x^{-n}(\Delta (\pi_1(S))) \cap \Delta(\pi_1(S))$.
This set is empty for
$\mid n \mid >> 1$,
so
$y =1$.
\enddemo
\demo{ Proposition 6.4}
(2nd criterion of virtual fibering). Let $M$ be a hyperbolic 3-manifold,
$\pi_1(S) \hookrightarrow \pi_1(M)$
an injective surface. Then $S$ is a fiber of a fibration for some finite covering $N$ of $M$ iff
\roster
\item"(i)" $\pi_1(S)$
is in infintely many subgroups of finite index and
\item"(ii)" for any
$x \in \pi_1(M)$
there exists
$N(x)$
such that for
$m \ge 1$
\endroster

\demo{Proof}
The "only if" part is obvious. On the other hand, by (i) $S$ is embedded in some covering $N$ of $M$. If $S$ is not a fiber of a fibration, then (ii) contradicts the Theorem A4.
\enddemo

\demo{6.5 Proof of the double coset theorem}
By the virtue of Sullivan' theorem we can find
$x \in \pi_1(M)$
such that its fixed points lie in different components of
$S^2 \setminus \Delta (\pi_1(S))$.
Say, its attracting point lie in $U$ and its repelling point lie in $L$ where 
$U,L$
are these components.
\enddemo

\demo{Claim}
For fixed positive $n$ and any 
$m >> n$,

$$x^m \notin \pi_1(S) x^n\pi_1(S)$$
\enddemo

\demo{Proof of the claim}
Suppose
$yx^n = x^mz$
with
$y,z \in \pi_1(S)$.
Let
$a \in U$
be the attractive fixed point of $x$. Since
$z(U) = U$,
for
$m >> 1$
the image
$x^mz(U)$
will be in a small neighbourhood of $a$, in particular its diameter in the hyperbolic metric of $U$ converges to $0$ as
$m \rightarrow \infty$.
On the other hand,
$\text{diam} \, (yx^n(U)) = \text{diam} (x^nU)$
since
$\pi_1(S)$
acts isometrically. Since $n$ is fixed,
$\text{diam} \, (yx^nU)$
does not converge to $0$.

Now one chooses inductively a sequence
$n_i$
such that
$x^{n_{\kappa}}$
does not belong to any of
$\pi_1(S) x^{n_i} \pi_1(S)$
with
$i \le k-1$.
The theorem is proved now.
\enddemo

\bf{7. Quasi-Fuchsian knots and analytic proof of Freedman-Cooper-Long theorem} 

\demo{7.1}
\enddemo
In this section we will apply the ideas developed earlier in connection with {\bf closed} surfaces in three-manifolds, to {\bf surfaces with boundary}. One of the outputs of this study will be an analytic proof of the Freedman-Long theorem discussed above
. The central notion in this section is that of quasi-Fuchsian knot. Recall first, that if $M$ is a three-manifold with boundary a torus $T, K$
a simple loop in $T, [K] =0$
in
$H_1(M,\Z)$
then there is a Seifert surface $S$ in $M$ with boundary at $K$, which is incompressible, that is,$\pi_1(S) \rightarrow \pi_1(M)$
is injective. Now suppose $M$ is hyperbolic.

\demo{7.2 Definition}
A knot $K$ is called quasi-Fuchsian if there is a Seifert surface $S$ such that the restriction of the uniformization representation
$\rho : \pi_1(M) \rightarrow PSL_2(\C)$
on
$\pi_1(S)$
is quasi-Fuchsian.
\enddemo

\demo{7.3 Remark }
The Thurston's proof of his hyperbolizatin theorem [Mo] implies that
{\bf all non-2-virtually fibered hyperbolic knots are quasi-Fuchsian} by Theorem  1.6.3.
\enddemo

\demo{7.4 Proof of the Theorem A3}
Let
$\ell = [K]$
in
$\pi_1(\Sigma \setminus N(\kappa)$,
and let
$m \in \pi_1(\Sigma \setminus N(\kappa))$
be a meridian. Let
$P \subset \pi_1(\Sigma \setminus N(K))$
be a group, generated by
$m, \ell$.
Then
$P \approx \Z \otimes \Z$.
Consider
$\rho (P) \subset PSL_2(\C)$.
We can conjugate
$\rho (P)$
into a group generated by
$\pmatrix 1 & 1 \\ 0 & 1 \endpmatrix$
and
$\pmatrix 1 & \tau \\ 0 & 1 \endpmatrix$
where
$Im \, \tau > 0$.

Now, let
$\Delta$
be a limit set for
$\rho (\pi_1(S))$.
This is a Jordan curve in
$\bar \C \simeq S^2$.
Since
$[\ell]$
is represented by the transformation
$z \mapsto z+1, \Delta$
is simply a periodic curve (Fig.1).

$$   $$
$$\Delta \qquad \qquad \qquad m \nearrow z \mapsto z+\tau \qquad \qquad \overset \ell \to {\underset z \mapsto z+1 \to \longrightarrow}$$
$$   $$
$$\text{Fig.1} \qquad$$
We will call the linear measure of the projection of
$\Delta$
on
$\C / \R \cdot 1$,
divided by
$\text{Im} \, \tau$,
a width of $K, w (K)$
for short. This is a topological invariant of $K$. Now we claim that for
$s \ge w(K)$,

$$m^s
\pi_1(S)m^{-s} \cap \pi_1(S) = \{\ell^n, n\in \Z\}$$
Suppose
$x \in m^s\pi_1(S)m^{-s}$.
Then the fixed points of $x$ lie on
$\Delta +s\cdot \tau$.
So if
$x \in m^s \pi_1(S)m^{-s} \cap \pi_1(S)$,
then the fixed points of $x$ lie on
$\Delta + s \cdot \tau \cap \Delta = \{\infty\}$,
so $x$ is parabolic therefore $x$ is conjugate to a power of
$\ell$,
and in fact $x$ is a power of
$\ell$
since
$Fix (x) = \{\infty\}$.
This proves the Theorem A3. It is elementary to deduce Theorem 3.2 from Theorem A6 [CL].
\enddemo

\bf{8. Freedman's surfaces give boundary groups} 

Let
$K \subset \Sigma$
be a hyperbolic knot which is not fibered. Let
$S^g$
be a minimal Seifert surface for $K$ and let $k$ be such that the Freedman's surface
$F = Fr(\pi_1(S), Ad \, m^{\kappa})$
is incompressible in the infinite cyclic covering of
$\Sigma \setminus K$
(by 3.2, it is enough to take
$k \ge 6g-2)$.
We may, and will, assume that the Freedman's surface is embedded in the cyclic covering
$M_{\kappa}$
of
$\Sigma \setminus K$.
Let
$\rho : \pi_1(\Sigma \setminus K) \rightarrow PSL_2(\C)$
be the uniformization representation. We wish to understand
$\rho \mid \pi_1(F)$.

\demo{8.1 Proposition}
$ \rho \mid \pi_1(F)$
is geometrically finite.
\enddemo

\demo{Proof}
Immediate by Bonahon's theorem [Bo].
\enddemo

\demo{8.2 Propositon}
$\rho \mid \pi_1(F)$
is a (Bers) boundary group.
\enddemo

\demo{Proof}
Let
$N = M_k$
and let
$\hat N_t$
be its Dehn fillings with slope
$(1,t)$.
Then, as in the proof of the Surgery Theorem A1, if
$\rho_t$
is the uniformization representation for
$\hat N_t$,

$$\rho_t\mid \pi_1(F) \rightarrow \rho \mid \pi_1(F)$$
Moreover, $F$ is incompressible in
$N_t[CL], [CGLS]$
therefore
$\rho_t\mid \pi_1(F)$
is quasi-Fuchsian, so
$\rho\mid \pi_1(F)$
is a boundary group.

For a deeper analysis of
$\rho\mid\pi_1(F)$
we notice that
$\hat N_t$
admits a natural decomposition

$$\hat N_t  = Y_0 \underset{\sigma} \to \cup Y_1$$
where
$Y_0, Y_1$
are defined as follows. Let $Y$ be $\Sigma -N(K)$ cut open along
$S$. Then $\partial Y=S_+ \cup  S_-$where $S_+ , S_-$ are canonically
homeomorphic to $S$. For any $d_0, d_1$, such that $d_0+d_1=t$  one
defines $Y_i=Y\cup Y\cup ....Y$($d_i$ times) always gluing $S_+$ of the previous copy of $Y$ with the
$S_-$ of the next. We assume $d_i\geq {6g-2}$, then $\partial Y_i$ is
incompressible in $Y_i$ by [CL]. The map 
$\sigma : \partial Y_0 \rightarrow \partial Y_1$
is a homeomorphism, described as follows:
$\partial Y_0$
is a double of
$S$,
so that there is a canonical orientation-reversing involution
$\lambda$,
fixing the middle curve $L$. Let
$\beta$
be the Dehn twist along $L$, then
$\sigma = \beta^t \circ \lambda$.
Let
$\theta_0 : T(\partial Y_0) \rightarrow T(\partial Y_0)$
and
$\theta_1 : T (\partial Y_1) \rightarrow T/\partial Y_1)$
be the skinning maps. Since
$Y_0$
is canonically homeomorphic to
$Y_1$
we can identify them. Then an element $z$ of
$T(\partial Y_0)$
which solves the gluing problem is a fixed point of
$\theta_0 \circ \sigma \theta_1  \sigma^{-1}$.
The situation is similar, but more complicated, than in [KT] where the authors looked at the boundary group defined by a sequence
$(z_0,\beta^tz_0)$
where
$z_0 \in T(\partial Y_0)$
is fixed and
$\beta$
as above. In [KT] the limit group has the limit set a proper subset of
$S^2$
whereas in our situation it seems that the limit set is all of
$S^2$
though I don't know how to prove this.
\enddemo

\bf{9. Digression to arithmetic topology: Freedman surfaces and generalized cyclotomic units} 

\demo{9.1}
The reader should be familiar with Appendix 21, The Language of Arithmetic Topology, in order to understand this section.

Let $K$ be a number field with
$p \nmid h(K)$.
Let
$L \supset K$
be a cyclic extension with
$\text{Gal} \, (L \mid K) =C_p$.
Then some prime ideal 
$\frak q \subset \Cal O (K)$
is totally ramified:
$\frak q = \frak f^p$
in $L$. Assume 
$\frak q$
is principal in $K$ (say
$h (K) = 0)$.
By the motivic spectral sequence,
$p \nmid h(L)$,
so $f$ is also principal. Let
$w \in \frak f$
be such that
$\frak f = (w)$.
Let
$\xi$
be a generator of
$\text{Gal} (L \mid K)$.
Then for any $m$,
$\xi^mw \cdot w^{-1}$
is a unit. This is true for any cyclic extension for which there is a totally ramified principal ideal 
$\frak q$
such that $\frak f$ is again principal. If
$K = \Q, L = \Q(p \sqrt 1)$,
so that
$\text{Gal} (L \mid K) = C_{p-1}$
then
$(p)$
is totally ramified:
$(p) = (1-\zeta)^{p-1}$
and one arrives to Kummer's cyclotomic unit
$1 +\zeta + \ldots + \zeta^{m-1}$

The topological counterpart is a 3-manifold
$\Sigma$
and its ramified covering
$\Sigma_m$.
The analogue of $w$ is the Seifert surface $S$, and
$\xi^mw \cdot w^{-1}$ 
is {\bf exactly} the Freedman's surfaces in
$\Sigma_m$.
\enddemo

\bf{10.Jaco's theorem and its corollaries}

Let $M$ be an acyclic closed 3-manifold, which does not contain a
surface  group. Then:
\demo{Proposition 10.2}
(Jaco) Any infinite index subgroup of
$\pi_1(N)$
is locally free and conversely.
\enddemo

\demo{Proof} 
Jaco [J].
\enddemo

\demo{Corollary 10.3}
There are maximal (normal) subgroups of infinite index in
$\pi_1(M)$,
i.e. there are infinite index (normal) subgroups
$\Gamma \subset \pi_1(N)$
such that any (normal)
$\Delta \supset \Gamma$
is of finite index.
\enddemo

\demo{Proof}
Immediate by Zorn lemma. The union of a chain of locally free subgroups is locally free.

That means {\bf there are infinite quotients $G$ of
$\pi_1(M)$
with all normal subgroups of finite index.}
\enddemo

\bf{11. $p$-cycles} 

\demo{11.1}
A $p$-cycle is a following 2-complex: take an oriented surface $S$ with boundary a circle
$\partial S$.
Take a product
$I \times W$
where $W$ is the graph in Fig.2.

\midspace{3cm}
$   $\newline
and identify
$\{0\} \times W$
with
$\{1\} \times W$
by rotation by
$\frac{2\pi} p$.
The resulting complex has a boundary 
$S^1$.
Glue it to
$\partial S$.
This is a $p$-cycle. For
$p=2$ a $p$-cycle is a nonorientable smooth surface.
\enddemo

\demo{11.2 Lemma}
For $p$ a prime and $X$ any space with
$H_2(X,\Z) =0$
any element in
$H_2(X,\F_p)$
can be represented by a map of a $p$-cycle.
\enddemo

\demo{Proof}
$H_2(X,\F_p) \approx Tor_1(H_1(X), \Z_p)$.
Let
$z \in \pi_1(X)$
be such that
$z^p = \prod^g_{i=1} [x_i,y_i]$.
so
$[z]$
represents  an element in
$Tor_1(H_1(x),\Z_p)$.
Map $W$ to $X$ sending the core loop to $z$ and send $S$ to $X$ by the map of generators and glue.
\enddemo

\demo{11.3 Lemma}
If $M$ is a closed 3-manifold and
$W \hookrightarrow M$
is an embedded $p$-cycle, then
$[Q] \in H_2(M,\F_p)$
is not zero.
\enddemo

\demo{Proof}
One constructs easily a loop in $M$ intersecting $Q$ transversaly at one point.
\enddemo

\demo{11.4 $p$-cycles in surgeries}
Let
$\Sigma$
be a homology sphere, $K$ a knot in
$\Sigma, S^g$
a Seifert surface of $K$ with boundary at
$\partial N(K) \approx T^2$.
Let
$\hat \Sigma (p,q)$
be a result of
$(p,q)$
Dehn surgery,
$(p,q) =1$.
\enddemo

\demo{Lemma 11.5}
$\Hat \Sigma (p,q)$
contains a $p$-cycle of genus $g$.
\enddemo

\demo{Proof}
This is imediate. The solid torus we glue contain a $W$-bundle over a circle with boundary parallel to $K$. Observe that
$H_1(\hat \Sigma (p,q)) \approx \Z/p \Z$
and the class of the $p$-cycle under consideration is dual to a generator of
$H^1(\hat \Sigma, \Z_p)$.
\enddemo

\demo{11.6 Cycles $mod \, p$ in ramified coverings}
Let
$\Sigma$ 
be a homology sphere and let $L$ be a link in
$\Sigma, L = K_1 \cup \ldots \cup K_s$.
Let $S$ be a Seifert surface of
$K_1$,
transversal to all
$K_i, i \ge 1$.
Let
$\Sigma_p$
be a $p$-ramified covering of
$\Sigma$.
Let
$\tilde S$
be a preimage of $S$ in
$\Sigma_p$.
Then
$\tilde S$
is a cycle
$mod \, p$
but not a $p$-cycle is in our sense, since the neighbourhood of the singular curve has a boundary which is connected. In particular,
$[\tilde S]$
may represent zero in
$H_2(\Sigma_p, \F_p)$.
However, if for some $i$, link
$(K_i,K_1) \ne 0 (mod \, p)$
then the class of
$[\tilde S]$
is not zero. This will be used later.
\enddemo

\demo{11.7 Definition}
A $p$-cycle
$Q \subset M$
is called dominating, if
$\pi_1(Q) \rightarrow \pi_1(M)$
has an image of infinite index in
$\pi_1(M)$.
\enddemo
\proclaim {Theorem 11.8}Let $M$ be a closed irreducible
three-manifold which has a non-dominating $p$-cycle.Then $ \pi_1(M)$
contains a fundamental group of a closed surface.
\endproclaim

We remark that if $p=2$, that is if $M$ contains a nonorientable
surface which is not dominating, then $M$ is Haken by a theorem of Rubinstein.
\demo{Proof}
Let
$Q \subset M$
be a non-dominating $p$-cycle. Then the image of
$\pi_1(Q)$
is of infinite index.
\enddemo

\demo{Claim}
The image $G$ of
$\pi_1(Q)$
in
$\pi_1(M)$
is not free.
\enddemo

\demo{Proof of the claim}
Suppose the opposite. Then the image of the fundamental class of $Q$ in
$H_2(\pi_1(M), \F_p) \approx H_2(M, \F_p)$
factors through
$H_2(G,\F_p) = 0$,
which is impossible.
\enddemo

The Theorem 11.8  follows now from Proposition 10.2.

\demo{Corollary 11.8}
Let $M$ be a closed manifold with
$H_1(M,\F_p) \ne 0$. Let
$z \in \pi_1(M)$
be such that
$z \notin [\pi_1(M), \pi_1(M)]$
but
$z^p \in [\pi_1(M), \pi_1(M)]$,
say
$z^p = \prod^g_{i=1} [x_i,y_i]$.
Suppose
$\{x_i,y_i,z\}$
generated an infinite index subgroup. Then
$\pi_1(M)$
contains a
$\pi_1$
-injective surface.
\enddemo

\demo{Proof}
Follows from above and 11.2.
\enddemo

\demo{Corollary 11.9}
Let
$\Sigma$
be a homology sphere $K$ a knot,
$S^g$
in Seifert surface and
$\hat \Sigma (p,q)$
the Dehn surgery.If no 
$2g+1$
elements generate a finite index subgroup of
$\hat \Sigma (p,q)$
then
$\hat \Sigma (p,q)$
contains an injective surface.
\enddemo

\demo{11.10 Remark}
If $K$ is not fibered, then for big $p$,
$p \ge 6g-2, \hat \Sigma (p,q)$
is virtually Haken by the Cooper-Long Theorem.
\enddemo

\demo{11.9 Remark}
Let
$\Sigma_{2p}$
be the 
$2p$-ramified
covering of
$\Sigma$.
If no
$2g$ elements generate
$\pi_1(\Sigma_{2p})$,
then obviously
$Y_0$
(see above) is not a handlebody, so a compressing of 
$\partial Y_0$
gives a closed incompressible surface in
$\Sigma_{\infty}$,
the infinite cyclic covering of
$\Sigma \setminus K$.
The theorem of Freedman 3.1 would follow then from the assertion that the least number of generators for
$\pi_1(\Sigma_{2p})$
goes to
$\infty$
as
$p \rightarrow \infty$,
if $K$ is not fibered. I do not know how to prove this.
\enddemo

\bf{12. Homology domination} 

\demo{12.1}
A $p$-cycle
$Q \subset M$
is $q$-homology dominating, if
$H_1(Q,\F_q) \rightarrow H_1(M,\F_q)$
is surjective.
\enddemo

\demo{12.2 Proposition}
Let
$Q \subset M$
be a $p$-cycle of genus $g$. If for some 
$q, b_1(M,\F_q) \ge g+2$,
then $M$ has a 
$\pi_1$
-injective surface.
\enddemo

\demo{Proof}
$Q$ is not dominating by [SW], $\qquad$. So $M$ has a 
$\pi_1$-injective surface by Theorem 11.6.
\enddemo

\proclaim{12.3 Theorem A5}
Let $M$ be a homology sphere with a link
$K = K_1 \cup \ldots \cup K_s$.
Let 
$N \rightarrow M$
be a $p$-covering, ramified along $L$. Let
$F_1$
be an oriented Seifert surface of
$K_1$,
transversal to
$K_i, i\ge 2$.
Let
$P = \sum^s_{i=2} \sharp (K_i \cap F_1)$.
Suppose for some
$i \ge 2$,
link
$(K_i,K_1) \ne 0(p)$.
If for some
$q \ne p$,

$$\dim H_1(N,\F_q) \ge 4-p(\chi(F)-2)+(p-1)P,$$
then $N$ is virtually Haken.
\endproclaim

\demo{Proof}
The full preimage $S$ of $F$ is a result of a gluing of a smooth (possibly disconnected) oriented surface $Q$ with $p$ boundary components and a product of a circle with the graph of Fig.2.

One sees immediately that $S$ is a cycle
$mod \, p$.
Since link
$(K_i,K_1) \ne 0$
(for some
$i, [S] \cap [K_i] \ne 0(p)$
in $N$, therefore,
$[S] \ne 0$
in
$H_2(N,\F_p)$.
It follows that the image of
$\pi_1(S)$
in
$\pi_1(N)$
is not free.

The Euler characteristic of $Q$ is bounded above by
$p\chi (F)-(p-1)P$
by Hurwitz formula. Since
$K_i \cap F \ne 0$
in $M$, the Galois group of the covering 
$N \rightarrow M$
may not permute the components of $Q$, so $Q$ is connected. It follows that
$\pi_1(Q)$
is generated by
$2-(\chi (F))p +(p-1)P+p$
elements.
\enddemo

\demo{Lemma}
Let
$X = Y\cup Z, A = Y\cap Z$
be CW-complexes so that
$\pi_0(Y) = \pi_0(Z) =0$.
Suppose
$\pi_1(Y)$
is generated by
$r_1$
elements,
$\pi_1(Z)$
is generated by
$r_2$
elements and
$\mid \pi_0(A) \mid =p$.
Then
$\pi_1(X)$
is generated by 
$r_1+r_2+(p-1)$
elements.
\enddemo

\demo{Proof}
Is elementary and left to the reader.
\enddemo

Now we see that
$\pi_1(S)$
is generated by
$2-(\chi (F))p+(p-1)P+p+1+p-1 =2-\chi (F)p+(p-1)P+2p$
elements. Now we may apply the theorem of Jaco-Baumslag-Shalen-Wagreich to prove the theorem.

\demo{12.4}
Now we consider the case
$p=q=2$.
We will prove the Theorem A4, stating that if $M$ posess a nonorientable surface which is not 2-homology dominating, then $M$ has virtually positive 
$b_1$.
We need to collect first some well-known facts about nonorientable surfaces in orientable 3-manifolds.
\enddemo

\demo{12.5}
Let $M$ be an orientable three-manifold and $S$ an embedded non-orientable surface in $M$. Recall that $S$ is diffeomorphic to a connected sum of
$2-\chi(M)$
copies of
$\R P^2$.
In particular
$w^2_1(S) = w_2(S) = \chi (M) (mod \, 2)$
since the Stiefel-Whitney numbers are bordism invariants. Let
$\nu$
be the normal bundle to $S$, then
$w_1(\nu) = w_1(S)$
since
$w_1(M) =0$.
Let
$\gamma$ be a simple curve in $S$, such that
$(w_1(S), [\gamma]) \ne 0$.
It follows that
$\nu \mid \gamma$
is nontrivial, therefore there is a section
$\delta$
of
$\nu \mid \gamma$
which is transversal to the zero section and intersects it in one point. That means that
$\delta$
is transversal to $S$ in $M$ and intersects it in pone point. In particular,
$[S] \ne 0$
in
$H_2(M,\F_2)$.
Moreover, for any simple curve
$\gamma$
in $S$ we have
$[\gamma] \cap [S] = (w_1(S), [\gamma])$,
where
$[\gamma]$
is understood to lie in
$H_1(M,\F_2)$
in the left hand side of the equation and in
$H_1(S,\F_2)$
in the right hand side.
\enddemo

Choose a section
$\Sigma$
of
$\nu$
over $S$ which is transversal to the zero section. The intersection
$\Sigma \cap S$
will represent a class in
$H_1(S,\F_2)$,
dual to
$w_1( \nu) = w_1(S)$.
It follows that
$S \cap S \cap S = w^2_1(S)$.
We resume this in the following proposition:

\demo{Proposition 12.5}
If
$S \subset M$
is as above then
$[S] \ne 0$
in
$H_2(M,\F_2)$.
Moreover, if
$\lambda \in H^1(M,\F_2)$
is a class dual to $S$, then
$\lambda^3 =\chi (S) (mod \, 2)$.
\enddemo

\demo{Corollary 12.6}
Suppose $M$ is a rational homology sphere and
$S \subset M$
is such that
$\chi (S)$
is odd. Then
$H_1(M)_{(2)} = \Z/2\Z \oplus V$
and the decomposition may be made orthogonal with respect to the linking form.
\enddemo

\demo{Proof}
By 4.1,
$\lambda^3 \ne 0$,
in particular
$0 \ne \lambda^2 = \beta (\lambda)$,
so that
$\lambda : H_1(M,\Z) \rightarrow \Z/2\Z$
does not factor though
$\Z/4\Z$.
It follows that
$H_1(M)_{(2)} = \Z/2\Z \oplus V$.
Moreover, if
$z \in H_1(M)_{(2)}$
is such that
$\lambda = (\cdot , z)$,
then
$\lambda^3 = (z,z) \ne 0$,
so taht the splitting may be made orthogonal.

Now consider an unramified double covering
$N \overset \pi \to \longrightarrow M$
corresponding to
$\lambda$.
Since we proved that
$\lambda \mid S = w_1(S) \ne 0$,
the preimage
$\pi^{-1}(S)$
is connected and oriented. If
$b_1(N) =0$,
then
$\pi^{-1}(S)$
should separate $N$ in two diffeomorphic manifolds
$D_1$
and
$D_2$
with boundary
$\pi^{-1}(S)$.
In other words, we have the following
\enddemo

\demo{Proposition 12.7}
Any three-manifold $M$ containing a non-orientable surface either has a double covering with positive
$b_1$,
or can be constructed in the following way: start with an oriented manifold $D$ with boundary
$\tilde S$, 
posessing on orientation-reversing free involution
$\zeta : \tilde S \rightarrow \tilde S$.
Then identify $x$ and $y$ if and only if
$x,y \in  \tilde S$
and
$\zeta x = y$.
\enddemo

\demo{Theorem 12.8}
Let
$p=2$.
If for some embedded non-orientable surface $S$, the natural map
$H_1(S,\Z_2)\rightarrow H_1(M,\Z_2)$
is not surjective, in particular, if 
$rank_2H_1(M)_{(2)} \ge b_1(S,\F_2)$,
then there is a double covering $N$ of $M$ with
$b_1(N) > 0$.
\enddemo

Here we will collect several facts from [Re 1] needed for a proof of the theorem.

\demo{12.9 Lemma}
([Re 1], 5.2). Let $N$ be a finite covering of $M$: Let
$\pi_*: H_1(N) \rightarrow H_1(M)$
be the natural map induced by
$\pi : N \rightarrow M$
and let
$t : H_1(M) \rightarrow H_1(N)$
be the transfer map. Let
$z \in H_1(M)_{tors}$
and
$w \in H_1(N)_{tors}$.
Then
$(\pi_*w,z)_M = (w,tz)_N$
where
$(\cdot , \cdot)$
is a linking form, valued in
$Q/\Z$.
\enddemo

\demo{12.11 Lemma}
([Re 1]). Let $p$ be a prime and let
$z \in H_1(M)_{(p)}$
be of order $p$. Let
$N \overset \pi \to \rightarrow M$
be a covering, corresponding to
$(\cdot , z)$.
Suppose
$b_1(N) =0$
then
$t(z) =0$,
and moreover
$\ker (t)$
is generated by $z$.
\enddemo

\demo{Proof of the theorem}
Choose
$\mu \in H^1(M,\Z_2)$,
different from
$\lambda$,
such that
$\mu \mid_S = \lambda \mid S$.
Let
$z,u \in H_1(M,\Z)_{(2)}$
be elements of order 2 such that
$(z,\cdot) = \lambda, (u,\cdot) = \mu$.
Let
$\pi : N \rightarrow M$
be a covering, corresponding to
$\mu$.
Suppose
$b_1(N) = 0$.
Let
$Q = \pi^{-1}(S)$.
\enddemo

\demo{12.12 Lemma}
Let
$w \in H_1(S,\Z)$
be unique nontrivial torsion element. Then the image of $w$ in
$H_1(M,\Z)$
is $z$.
\enddemo

\demo{Proof}
Let
$\lambda \subset M$
be a closed curve, transversal to $S$. There exists a decomposition
$S = S_0 \cup \Sigma_i$
such that
\roster
\item"(a)" $S_0$
is a sphere with
$g = genus \, (S)$
discs removed. We denote
$C_i$
the boundaries of these discs
\item"(b)" $\Sigma_i$
is a M\"obius band with basis 
$d_i$
and boundary 
$C_i$
\item"(c)" $w$ is represented by
$\Sigma d_i$
\item"(d)" all intersection points of
$\gamma$ 
and $S$ lie in
$S_0$.
\endroster
Let $v$ be a class of
$\Sigma d_i$
in
$H_1(M,\Z)$.
Then 
$2v$
bounds an oriented singular surface $B$ (a map at sphere with $g$ holes in $M$) with image $S$, and
$([\gamma], v)_M = 1/2 \sharp (\gamma \cap S_0)(mod \, \Z)$,
where
$\sharp (\gamma \cap S_o)$
is counted with orientation. Up to the identification
$\Z \cdot \frac 1 2(mod \, \Z) \approx \Z_2$,
this is
$[\gamma] \cap [S_0] = \lambda (\gamma)$.
\hfill q.e.d.
\enddemo

Now we have a commutative diagram

$$\CD H_1(Q) @>>> H_1(N) \\
@At AA @AAt A \\
H_1(S) @>>> H_1(M) \endCD$$
Since $Q$ is oriented by 14.5,
$H_1(Q)_{tors} = 0$,
so
$t(w) =0$.
It follows that
$t(z) =0$.
On the other hand, by Lemma 14.11,
$\ker (t)$
is generated by
$w \ne z$,a contraction. So
$b_1(N) > 0$.

\bf{13. Tight knots}

\demo{13.1 Definition}
A homologically trivial knot $K$ in a rational homology sphere
$\Sigma$
is called tight, if there is a Seifert surface $S$ for $K$, such that
$\pi_1(S) \rightarrow \pi_1(\Sigma)$
is injective.
\enddemo

\demo{13.2 Examples}
1. Let
$\Gamma \subset \pi_1(\Sigma)$
be any finitely generated free group (recall that any 3-manifold with rich fundamental group [Re 1] has free subgroups of any rank in its fundamental group). Realize
$\Gamma$
by a handlebody
$P \subset \Sigma$.
Let
$K \subset \partial P$
be any essential loop such that
$[P] =0$
in
$H_1(P,\Z)$.
Then there is an embedded incompressible surface $S$ in $P$ with boundary in a neighbourhood of $K$. So
$K = \partial S$
is tight.

2. If
$K \subset \Sigma$
is not fibred, then $K$ becomes tight in some ramified covering
$\Sigma_m$
of
$\Sigma$,
by [CL]
.

3. Let
$\Sigma$
be hyperbolic and such that no 2 elements generate a finite index subgroup in
$\pi_1(\Sigma]$.
Let
$K \subset \Sigma$
be a genus 1 knot,
$[K] \ne 1$
in
$\pi_1(M)$.
Then either
$\Sigma$
contains a
$\pi_1$
-injective closed surface, or $K$ is tight. Indeed, if $S$ is a punctured torus with
$\partial S = K$.
Suppose
$\psi : \pi_1(S) \rightarrow \pi_1(\Sigma)$
is not injective. It its image is free of rank 2, then
$\psi$
is an isomorphism [LS], a contradiction. If
$\psi (\pi_1(S))$
is free of rank one, then
$[K]=1$
in
$\pi_1(\Sigma)$.
If
$\psi(\pi_1(S))$
is not free, then
$\Sigma$
contains a 
$\pi_1$-injective surface.

4.A recent theorem of Boyer-Culler-Shalen -Zhang  [BCSZ] asserts that if $K
\in  \Sigma$ is not fibered, then for $m$ big enough , $K$ bevomes
tight 
in a $(m,n)$-surgery of  $\Sigma$.
\enddemo

\demo{Proof of the Theorem A6}
Consider the Freedman surface
$S_0 \cup \tau^{\ell}S_0, \ell\ge 1$.
We will use the technique of {\bf collapse} which consists in the
following lemma of Shapiro and Sonn [SS]  (a stonger theorem goes back
to Zieschang, see [Re]).
\enddemo

\demo{Lemma}
If
$\pi_1(S^g ) \overset \psi \to \longrightarrow \pi_1(F)$
is a surjective way of a fundamental group of a closed surface on a free group, then
$rank \, F \le g$.
\enddemo

\demo{Proof}
The image of
$H^1(F,\Q)$
in
$H^2(S,\Q)$
is Lagrangian with respect to the natural symplectic structure in
$H^2(S,\Q)$,
because
$H_2(F) =0$,
and so for any two elements in
$H^1(F,\Q)$,
say 
$u,v$,
$(\psi^*u \cup \psi^*v) ([S]) = (u \cup v) \psi_* [S] =0$.

Coming back to the proof of the theorem, we have the possibilities
\roster
\item"(a)" $\psi_1(S_0 \cup \tau^{\ell} S_0) \rightarrow \pi_1(\Sigma_m)$
has the image of finite index. Projecting to
$\Sigma$,
we see that
$\pi_1(S) \rightarrow \pi_1(\Sigma_m)$
has the image of finite index, a contradiction to tightness.
\item"(b)" the image of
$\pi_1(S_0 \cup \tau^{\ell} S_0) \rightarrow \pi_1(\Sigma_m)$
is free. By the lemma above, this is a free group $F$ with
$\le 2g$
generators. Now, $F$ projects surjectively on
$\pi_1(S) = F_{2g}$,
do $F$ has exactly $2g$ generators and the projection is an isomorphism [LS]. But the restrictions of the projection to
$\pi_1(S_0)$
and
$\pi_1(\tau^{\ell}S_0)$
is obviously an isomorphism as well. It follows that
$\tau^{\ell}_*\mid\pi_1(S_0)$
is identity. However, by a theorem of Conner [Co], for an action of a finite
cyclic group on a $K(\pi , 1)$  three-manifold space, $\tau^{\ell}_*$ cannot have fixed
elements in the fundamental group, except powers of the element,
represented by $K$. So this is impossible.
\item"(c)" the image of
$\pi_1(S_0 \cup \tau^{\ell} S_0) \rightarrow \pi_1(\Sigma_m)$
is not free and of infinite index. Then
$\Sigma_m$
has a 
$\pi_1$
-injective surface by 12.2. Now, we start to compress $ S_0 \cup
\tau^{\ell} S_0$. Then either we arrive to an incompresible surface, or
the map $\pi_1(S_0 \cup \tau^{\ell} S_0)\to \pi_1(M)$ is factored
through a free group $F_{2g}$ with $2g$ generators. Say, $f$ is a map
from $\pi_1(S_0 \cup \tau^{\ell} S_0)$ to $F_{2g}$ and $g$ is a map
from $F_{2g}$ to $\pi_1(\Sigma_ m)$. Arguing exactly as above, we find that $g$
must be an isomorphism, and $\tau^{\ell}_*\mid\pi_1(S_0)$ is identity, which
we showed to be impossible. So $\Sigma_ m$ is Haken.
\endroster
\enddemo

\demo{A7. Corollary}.Let
$K \subset \Sigma$
be any knot which is not fibered.Let $\hat \Sigma(m,n)$ be a
result of the Dehn surgery.Let $\hat \Sigma(m,n,p) $  be a
$p$-ramified covering of  $\hat \Sigma(m,n)$.Then for $m >> 1$ and
$p\geq 2$, $\hat \Sigma(m,n,p) $ is Haken or reducible.
\enddemo

\bf{14. Appendix: The Language of arithmetic topology} 

The circle of ideas which form the basis of what follows is a joint invention of several generation of mathematicians. The fundamental idea - that a closed three-manifold "is" a number field is a direct development of the classical analogy between arithme
tics and algebraic geometry. In various form this idea has been suggested by Mazur, Manin and others starting from 60-ies. Its present formulation is closed to that contained in a lecture of Kapranov and myself in MPI, in the summer of 1996.

\demo{14.1.1 Foundations}
A closed {\bf three-manifold} is an analogue of a {\bf number field}.
\roster
\item"2." A {\bf three-sphere}
$S^3$
is an analogue of $\Q$.
\item"3." A {\bf ramified covering}
$M \rightarrow N$
is an analogue of an {\bf extension}
$K \subset L$.
\item"4." A {\bf Galois ramified covering}
$M \overset G \to \longrightarrow N$
is an analogue of a {\bf normal extension}
$K \subset L, G = Gal (L \mid K)$
\item"5." An {\bf element}
$w \in \Cal O(K)$
is an analogue of an {\bf embedded surface with boundary}
$S \subset M$
\item"6." A {\bf knot}
$K \subset M$
is an analogue of a {\bf prime ideal}
$\frak p \subset \Cal O(K)$.
\item"7." A map
$w \mapsto (w)$
is an analogue of a correspondence
$S \mapsto \partial S$
\item"8." A {\bf closed embedded surface} 
$F \subset M$
is an analogue of a {\bf unit}
$\varepsilon \subset U(K)$
\item"9." {\bf Normalization of ramified coverings, ramified loci}, etc. correspond in the obvious way.
\item"10." The Galois group of a maximal unramified extension corresponds to the fundamental group.
\endroster
\enddemo

\demo{14.2 One step deeper}
\enddemo

\demo{14.2.1}
An ideal class group
$I(K)$
is an analogue of
$H^{tors}_1(M,\Z)$
\enddemo

\demo{14.2.2}
For a Galois ramified covering
$M \rightarrow N$
(extension 
$L \mid k)$
the Galois modules
$I(K), (H^{tors}_1(M,\Z))$
correspond.
\enddemo

\demo{14.2.3}
The group of units mod torsion
$U(K)/tors$
corresponds to
$H_1(M,\Z)/tors$.
\enddemo

\demo{14.2.4}
An embedding
$K \rightarrow \R$
or
$\C$
correspond to a harmonic 1-form
$w \in \Omega^{\prime}(M,\R)$
with respect to a hyperbolic metric of $M$ (if $M$ is hyperbolic).
\enddemo

\demo{14.2.5}
Freedman's surfaces correspond to (generalized) cyclotomic units, see
section 9.
\enddemo

\demo{14.2.6}
$K$-theory
$K_i(K)$
correspond to
$\pi_0(Diff (\undersetbrace i \to {M \times\ldots \times M}))$
\enddemo

\demo{14.3 Two steps deeper}
\enddemo

\demo{14.3.1}
A theorem of Gauss, stating that for $K$ a quadratic field,
$\dim_{\F_2}\Z/2\Z \otimes I(K)$
equals number of primes in
$D(K)$ minus one,
corresponds to Theorem 17.1.2.
\enddemo

\demo{14.3.2}
A theorem of Redei [Red] stating that for
$D(K) =p_1p_2$
($K$ quadratic),
$I(K)_{(2)}$
in
$\Z/2\Z$
or
$\Z/2^n\Z, n \ge 2$
depending on
$\pmatrix p_1 \\ p_2 \endpmatrix$
corresponds to Theorem A9.
\enddemo

\demo{14.3.3}
Further development see in [Re6].
\enddemo

\demo{14.3.4}
A solution by $\check S$afarevi$\check c$ of the class tower problem
of Furtw\"angler corresponds to the part (a) of the Theorem 10.1  of [Re 1], due essentially to Turaev [T].
\enddemo

\demo{14.3.5}
A Theorem 10.2 of [Re 1] states that if $M$ does not have virtually positive $b_1$
and
$b_1(M,\F_p) \ge 4$,
then the pro-p completion of
$\pi_1(M)$
is a pro-p Poincar$\acute e$ duality group.
\enddemo

\bf{15. Hyperbolization of Heegard splittings and bounded width in the mapping class group} 

\demo{15.1}
Thurston conjectured in [Th1] that "big" Heegard splittings are hyperbolic. This indeed seems to be the case, thou our approach is not that proposed by Thurston in [Th1].

Recall that the mapping class group
$Map_g$
is generated by Dehn twists (by Dehn and Likorish). In many ways 
$Map_g$
behaves as a lattice in a semisimple Lie group. For arithmetic lattices
$\Gamma$
of rank
$\ge 2$
it is porved in many cases, that they have a bounded width in the following sense: there is a finite set 
$\{\lambda_i\}^N_1 \in \Gamma$
such that any element is of the form
$\prod^N_{i=1} \lambda^{k_i}_i$.
\enddemo

\demo{Conjecture}
(Bounded width in
$\text{Map}_g$).
There are finitely many Dehn twists
$\lambda_i$
in
$\text{Map}_g$
such that any double coset in

$$\text{Map}^0_g \setminus \text{Map}_g / \text{Map}^0_g$$
has a representative of the form
$\prod^N_{i=1} \lambda^{m_i}_{j_i}$,
for some fixed
$N \in \N$.
\enddemo

Here one fixes a handlebody $Y$ with
$\partial Y = S^g$
and denotes
$Map^0_g$
the image in
$\pi_0(Diff (S))$
of
$\pi_0(Diff (Y))$.

Let
$\{L_i\}^N_{i=1}$
be a family of simple closed essential curves in $S$. Consider
$Q = S \times [0,N]$
and let
$L =  \cup L^{(i)}$
be a link in $Q$ where
$L^{(i)}$
is a copy of
$\{L_i\}$
put on
$S \times \{i\}$.
Let
$Y_0 \cup Y_1$
be a Heegard splitting of
$S^3$
(they are all standard). Thickening
$\partial Y_0$
we represent
$S^3$
as
$Y_0 \cup S \times [0,N] \cup Y_1$.

\demo{Proposition 15.3}
For a set of integers
$\{k_i\}$,
a result of the surgery on the link $L$ with parameters
$(1,k_i)$
is diffeomorphic to the Heegard splitting with the gluing
$\prod {\beta_i}^{k_i}$
where
$\beta (L)$
is the Dehn twist around $L$.
\enddemo

\demo{Proof}
An exercise left to the reader.
\enddemo

Now, for $k_i$ all big enough, the Thurston conjecture above  follows
from the Thurston surgery theorem.

\bf{Appendix 16. Euler class and free generation} 

%$   $\newline
%``What can you say in your favour?'' \newline
%``You see ...'' \newline
%``Enough. Shoot'm. Next.'' [S].

%$   $\newline

This appendix consists of two parts. In the first auxilliary part, we deal with sets with cyclic order. For any such set
$\Cal O$
we introduce following [BG] a cocycle
$\ell : G \times G \rightarrow \Z$
valued in
$\{0, \pm 1\} \subset \Z$
on the group $G$ of automorphisms of
$\Cal O$.
The cohomology class of
$\ell$
in
$H^2(G,\Z)$
will be called the Euler class. If $K$ is an ordered field, then the projective line
$\P^1(K)$
has a cyclic order and
$PSL_2(K)$
acts order-preserving on
$\P^1(K)$, so that we get both the cocycle 
$\ell$
and the Euler class in
$H^2(PSL_2(K),\Z)$.
If
$K = \R$,
then our Euler class coincides with the usual Euler class on
$PSL_2(\R)$.

In view of our extension of the Euler class to all ordered fields, the following two problems arise.

\demo{{\bf Problem 1}}
Let
$\rho : \pi_1(S) \rightarrow PSL_2(K)$
be a homomorphism of the fundamental group of a closed oriented surface. Is it true, that
$|(\rho^*[\ell], [S])| \le 2g-2$?
\enddemo

\demo{{\bf Problem 2}}
Suppose
$|(\rho^*[\ell], [S])| = 2g-2$.
Is it true that
$\rho$
is injective?
\enddemo

For
$\K = \R$
the theorems of Milnor [M] and Goldman [Go2], answer these problems positively.

In the main second part of this paper, we apply the cocycle
$\ell$
and the ideas from the theory of Hamiltonian systems on the Teichm\"uller space to the following classical problem:

\demo{{\bf Problem 3}}
When $n$ matrices in
$SL_2(\R)$
generate a free discrete group?
\enddemo

For
$n=2$
this problem has been treated in many papers see [Gi]. An effective solution is given in [Gi]. The analogous problem for
$SL_2(\C)$
also attracted a lot of attention especially since J\"orgensen paper [J]. For
$n>2$
however, the problem becomes much harder. We will give a simple {\it sufficient} condition for $n$ hyperbolic elements in
$SL_2(\R)$
to generate a free discrete group. This condition is open, that is, satisfied on an open domain in
$(SL_2(\R))^n$.
Here is our main result.

Let $n$ be even and let
$a_i,b_i, 1\le i \le n$
be in
$SL_2(\R)$.
Suppose
$h = \prod^n_{i=1} [a_i,b_i]$
is hyperbolic. Consider the eigenvectors
$x_1, x_2$
of $h$ and a matrix $r$ which takes the form
$\pmatrix 1 & 0 \\ 0 & -1 \endpmatrix$
in the basis
$(x_1, x_2)$.
Put
$a_i = rb_{n+1-i}r^{-1}, b_i = ra_{n+1-i}r^{-1}$
for
$n+1 \le i \le 2n$.
Let
$I_j = \prod^j_{i=1} [a_i,b_i] (j \le 2n)$.

\proclaim{16.0. Theorem}
Let
$f(a,b) = \frac 1 {\pi} \sum^{2n}_{j=1} \ell(I_{j-1},a_j) + \ell (I_{j-1}a_j,b_j)- \ell (I_{j-1}a_jb_j a^{-1}_j,a_j) - \ell (I_j,b_j)$
Then
\roster
\item"(a)" $f(a,b)$
is an integer and
$|f(a,b)| \le 2n-1$
\item"(b)" if
$|f(a,b)| = 2n-1$,
then
$\{a_i,b_i\}$
generate a free hyperbolic group in
$SL_2(\R)$.
\endroster
\endproclaim

\demo{\bf{16.1. Cyclically ordered sets, ordered fields and the Euler class}} 
\enddemo

\demo{16.1.1}
A cyclically ordered set
$\Cal O$
is a set with a subset
$\Omega$
in
$\Cal O \times \Cal O \times \Cal O$,
satisfying the following conditions:
\roster
\item"(i)" if
$(x,y,z) \in \Cal O$
then
$x,y,z$
are all different
\item"(ii)" if$\sigma$
is a permutation of
$(x,y,z)$
and
$(x,y,z) \in \Cal O$
then
$\{\sigma (x), \sigma (y), \sigma (z)\} \in \Cal O$
if and only if $\sigma$ is even
\item"(iii)" if $z$ is fixed then the relation
$x < y \Leftrightarrow (x,y,z) \in \Cal O$
is a linear order.
\endroster
\enddemo

\demo{16.1.2 Example} Let $K$ be an ordered field and let
$\P^1(K)$
be a projective line over $K$. We can think of
$\P^1 (K)$
as
$K \cup\{\infty\}$.
The cyclic order in
$\P^(K)$
is defined by a condition that the induced order in $K$ is standard. The group
$PSL_2(K)$
acts on
$\P^1(K)$
preserving the cyclic order.
\enddemo

\demo{16.1.3}
Define a function
$\psi : \Cal O \times \Cal O \times \Cal O \rightarrow \{0, \pm 1\}$
in a following way:
\roster
\item"(i)" if any of
$(x,y,z)$
are equal, then
$\psi (x,y,z) =0$
\item"(ii)"
$\psi$
is odd under permutation of
$(x,y,z)$
\item"(iii)" if
$(x,y,z) \in \Omega$,
then
$\psi(x,y,z) =1$.
\endroster
\enddemo

\demo{16.1.4} Now let $G$ be a group, acting in order preserving way on
$\Cal O$.
Fix any element
$p \in \Cal O$
and define a function
$\ell : G \times G \rightarrow \{0, \pm 1\}$
as
$\ell(g_1,g_2) = \psi (p,g_2p,g_1g_2p)$.
\enddemo

\proclaim{Lemma (16.4)}
$\ell$
is an integer cocycle on $G$.
\endproclaim

\demo{Proof}
is a direct computation and left to the reader.
\enddemo

\demo{Definition}
The cohomology class of
$\ell$
in
$H^2(G)$
(which does not depend on $p$) is called the Euler class. In particular, for an ordered field $K$ one gets the Euler class in
$H^2(PSL_2(K), \Z)$.
\enddemo

\proclaim{16.1.5 Comparison theorem ([BG])}
For
$K = \R$,
the class of
$\ell$
in
$H^2(PSL_2(\R))$
is the usual Euler class of associated
$S^1$-bundle over
$BPSL_2^{\delta}(\R)$.
\endproclaim

\demo{Proof}
Consider the action of
$PSL_2(\R)$
on
$\Cal H^2$.
For any
$p \in \Cal H^2$
the class
$\ell_p(g_1,g_2) = \text{Area}\, (p,g_2p,g_1g_2p)$
represents the Euler class $e$ [Gu]. Here
$\text{Area} \, (p,q,r)$
is the area of oriented geodesic triangle with vertices in
$p,q,r$.
Now
$[\ell_p] \in H^2(PSL_2(\R))$
does not depend on $p$ and all cocycles 
$\ell_p$
are uniformly bounded. For
$p_0 \in \partial \Cal H^2$
and
$p \rightarrow p_0$,
we will have
$\ell^{\infty}$-convergence of cocycles
$\ell_{p_i} \rightarrow \ell_{p_0}$.
Moreover the area of ideal triangle
$(p,q,r)$
is
$\pi \cdot \psi(p,q,r)$.
Any homology class in
$H_2(PSL_2(\R))$
is represented by a map of a surface group
$\pi_1(S) \overset{\alpha} \to \longrightarrow PSL_2(\R)$.
It follows that
$(\ell_{p_0}, \alpha_*[S]) = \lim(\ell_{p_i}, \alpha_*[S]) = (e,\alpha_*[S])$.
This completes the proof.
\enddemo

\demo{\bf {16.2. Discrete Goldman twist}}
We will work with the representation variety
$\Cal M = \text{Hom} \, (\pi_1(S), SL_2(\R))/SL_2(\R)$,
where $S$ is an oriented closed surface of genus $g$. It has
$1+2(g-1)$
connected components, indexed by the value of the Euler class [H]. Every such component, say
$\Cal M_e$,
is a symplectic manifold, nonsingular if
$e \ne 0$.
For any
$\gamma$
a conjugacy class in
$\pi_1(M)$,
there is a natural Hamiltonian 
$Tr_{\gamma} : \Cal M \rightarrow \R$,
and the corresponding Hamiltonian flow has been identified by Goldman. If
$\gamma$
can be represented by a simple separating curve, this flow can be described as follows. Write a presentation of
$\pi_1(S)$
in the following form:

$$[x_1,x_2] \ldots [x_{2\kappa -1},x_{2\kappa}] = [\gamma] = [x_{2\kappa+1}, x_{2\kappa+2}] \ldots [x_{2g-1},x_{2g}]$$
Next, write
$[\gamma] = \exp A$
for some
$A \in sl_2(\R)$.
Then put

$$\aligned \bar x^t_i &= \bar x_i, i \le 2\kappa \\
\bar x^t_i &= \exp(-tA) \bar x_i \exp tA, i \ge 2\kappa; \endaligned$$
this is the flow of
$Tr_{\gamma}$.
In particular,
$f_t: \{\bar x_i\} \rightarrow \{\bar x^t_i\}$
is a symplectomorphism of
$\Cal M$.
Here
$\bar x_i$
stands for the representation matrix of
$x_i$.

For different
$\gamma$,
the Hamiltonians
$Tr_{\gamma}$
yield nice commutation relations, discovered by Wolpert [W] and put in a more ``representation variety language'' by Goldman [Go1]. In fact, the integer group ring of
$\pi_1(M)$
becomes a Lie ring with Goldman's bracket. One may wonder what kind of group object correspond to it.

Whatever this eventual ``Kac-Moodi-Goldman'' group may be, we will introduce now some elements from the ``other connected components'' of it. These are defined for
$\Cal M_{\pm (g-1)}$,
which is naturally symplectomorphic to the Teichm\"uller space by a theorem of Goldman (see various proofs in [Go2]. In this case, all representation matrices are hyperbolic.

For
$\gamma$
as above, let $r$ be a unique matrix (up to sign), commuting with
$[\bar {\gamma}]$
with eigenvalues
$+1$
and
$-1$.
Put

$$\aligned f(\bar x_i) &= x_i, i \le 2\kappa \\
f(\bar x_i) &= r^{-1} \bar x_i r, i \ge 2\kappa \endaligned$$
This is a symplectic diffeomorphism of
$\Cal M_{\pm(g-1)}$.
There is a particularly nice description of the map $f$ if one views
$\Cal M_{g-1}$
as Teichm\"uller space. Namely, realise a point in
$\Cal M_{g-1}$
as a hyperbolic metric on $S$ and find a geodesic, representing
$\gamma$.
We assume that the marked point lies on
$\gamma$.
Cut $S$ into two pieces along
$\gamma$
and glue again by a reflection, which fixes a marked point. Then the new hyperbolic structure is the image of $f$.
\enddemo

\demo{\bf{16.3. Euler class}}
For a representation
$x_i \rightarrow \bar x_i$
of
$\pi_1(S)$
in
$SL_2(\R)$
the Euler number is an integer between
$-(g-1)$
and
$(g-1)$
by Milnor [M]. As mentioned above, all representation with the maximal Euler number are discrete faithful hyperbolic by Goldman's theorem. One can introduce a universal Euler class $e$ in
$H^2(SL^{\delta}_2(\R),\R)$
as the image of a generator in continuous cohomology
$H^2_{cont}(SL_2(\R),\R)$.
A representation
$x_i \rightarrow \bar x_i$
as above defines a homology class in
$H_2(SL^{\delta}_2(\R),\Z)$,
the image of the generator of
$H_2(\pi_1(S),\Z) \approx \Z$,
and the Euler number is just given by the evaluation of $e$ on this class. Now, the generator of
$H_2(\pi_1(S),\Z)$
can be realized by an explicit cycle in the standard complex [Br], that is,
$\sum^{2g}_{j=1} (I_{j-1}|x_j) + (I_{j-1}x_j|y_j) -(I_{j-1}x_j y_jx^{-1}_j|y_j) - (I_j|y_j)$
where
$I_j = [x_1,y_1] \ldots [x_j, y_j]$.
On the other hand, the universal Euler class may be realized by a cocycle
$A,B \mapsto \ell (A,B)$
as defined in Section 1. So the Euler number will be

$$\sum^y_{j=1} \ell (\bar I_{j-1}, \bar x_j) + \ell (\bar I_{j-1}, \bar x_j, \bar y_i) - \ell (\bar I_{j-1} \bar x_j \bar y_j \bar x^{-1}_j, \bar y_j) - \ell (\bar I_j, \bar y_j),$$
where
$\bar I_j$
is defined in an obvious way.
\enddemo

\demo{\bf{16.3. Proof of the Main Theorem}}
Consider a closed surface $S$ of genus 
$2n$. 
A map
$x_i \rightarrow a_i, y_i \rightarrow b_i, 1 \le i \le 2n$
defines a homomorphism from
$\pi_1(S)$
to
$SL_2(\R)$.
Indeed,
$\prod^n_{i=1} [a_i,b_i] \cdot \prod^{2n}_{i=n+1} [a_i,b_i] = h \cdot r^{-1} h^{-1}r =1$.
Next, the Euler number of this representation is computed as above, so
$f(a,b)$
is always an integer. Moreover, if
$f(a,b) = 4n-2$,
then the representation above is discrete and faithful.
\hfill Q.E.D.
\enddemo

\demo{\bf{$SU(1,n)$-case, I}}
Consider a standard action of
$SU(1,n)$
on the unit ball
$B \subset \C^n$
with the complex hyperbolic metric. Let
$\omega$
be the K\"ahler form of $B$. Fix a point
$\infty$
in the sphere at infinity and consider a function
$\varphi(A,B) = \psi(\infty, A(\infty), AB(\infty))$,
where
$\psi(x,y,z)$
is an integral of
$\omega$
over any surface, spanning the geodesic triangle with vertices
$x,y,z$.
Let
$a_i,b_i, \, 1\le i \le g$
be matrices in
$SU(1,n)$.
Let
$h = \prod^g_{i=1} [a_i, b_i]$
and suppose $h$ has a nonisotropic eigenvector. Then there exists a reflection $r$, commuting with $h$. Define
$a_i, b_i, \, i \ge g+1$
as in Introduction.
\enddemo

\proclaim{16.3.1 Theorem}
Define
$f(a,b)$
by the formula in the Main Theorem. Then
\roster
\item"(a)" $f(a,b)$
is an integer and
$|f(a,b)| \le 2g-1$
\item"(b)" if
$f(a,b) = 2g-1$,
then
$\{a_i,b_i\}$
generate a discrete group in
$SU(1,n)$.
\endroster
\endproclaim

\demo{Proof}
Same as above with Toledo's theorem [To] instead of Goldman's.
\enddemo

\demo{\bf {16.4. Bounded cohomology}}
The bounded cohomology theory is an invention of Mikhail Gromov. The idea is as follows: in the standard complex, computing the real group cohomology we look only at bounded cochains, that is, bounded functions
$f : \undersetbrace i \to {G \times G \times \ldots \times G} \rightarrow \R$.
The resulted cohomology spaces are called
$H^i_b(G,\R)$.
There is a canonical homomorphism
$H^i_b(G,\R) \rightarrow H^i(G,\R)$.
\enddemo

\demo{16.4.1 Example}
Let $M$ be a symmetric space of negative curvature with isometry group $G$. Let
$\omega$
be any $G$-invariant $i$-form on $M$. Then one gets a (Borel) class
$Bor \, (\omega) \in H^i_{cont}(G,\R)$
(see [Re] for example), which may be represented by bounded cocycle (Gromov [Re1]). The Euler class in
$SU(1,n), \, n \ge 1$
is a further specialization.
\enddemo

\demo{16.4.2 Second bounded cohomology and combinatorial group theory}
Consider a kernel of the map
$H^2_b(G,\R) \rightarrow H^2(G,\R)$.
It gives rise to a function
$f:G \rightarrow \R$
satisfying
$|f(x,y) -f(x)-f(y)| \le C$.
Moreover, this function may be chosen a class function, that is,
$f(xyx^{-1}) = f(y)$
and such that
$f(x^n) = nf(x)$
[BG].
\enddemo

Next, for an element
$z \in G^{\prime} = [G,G]$
a {\it genus norm} is the smallest integer $g$ such that $z$ is a product of $g$ generators. A theorem of Culler [C] states:

\proclaim{Theorem 18.4.2 (Culler)}
If $G$ is a f.g. free group, then for any
$z \in G^{\prime}$,
$\parallel z^n\parallel_{genus} \ge const \, \cdot n$.
\endproclaim

The following result has been proved by Gromov [Gr2] and author [Re3].

\proclaim{Theorem 16.4.3}
Let $G$ be geometrically hyperbolic, that is a fundamental group of a manifold of pinched negative curvature with
$i(x) \underset{x\rightarrow \infty} \to \longrightarrow \infty$
(e.g. compact). Then the conclusion of the Theorem 4.2 holds.
\endproclaim

Now, we have

\proclaim{Proposition 16.4.4}
Let
$f:G \rightarrow \R$
be as above. If
$f(z) \ne 0$,
then
$\parallel z^n \parallel_{genus} \ge const \, \cdot n$.
\endproclaim

\demo{Proof}
Let
$z^n = \prod^g_{i=1} [x_i,y_i]$.
Then

$$\gather |n \cdot f(z)| = |f(z^n)| = |f(\prod^g_{i=1} [x_i,y_i])| = |\sum^g_{i=1} f([x_i,y_i])| + \\
+ C \cdot g \le 3c \cdot g + C \cdot g = 4C \cdot g, \quad \text{so} \quad g \ge \frac{n|f(z)|} {4C}, \endgather$$
\hfill Q.E.D.
\enddemo

\demo{\bf{Genus norm in lattices in $SU(1,n)$}} \enddemo

\proclaim{Theorem (16.4.5)}
Let
$\Gamma \subset SU(1,n)$
be a lattice with
$H_2(\Gamma,\R) =0$. 
There exists a nonzero function
$f : \Gamma \rightarrow \R$
as above such that if
$f(z) \ne 0$,
then
$\parallel z^n \parallel_{genus} \ge const \, \cdot n$.
\endproclaim

\demo{Remarks}
1. Observe that if
$f(z) \ne 0$,
then for any
$y \in G$
and
$\kappa$
big enough,
$f(z^{\kappa}y) \ne 0$,
so there are ``many'' $z$ for which the Theorem 6.1 applies. \newline
2. If
$\Gamma$
is cocompact or if
$\underset{\gamma \rightarrow \infty} \to {|Tr \, \gamma|} \rightarrow \infty$
in
$\Gamma$
then the conclusion of the Theorem follows from 16.4.3.
\enddemo

The proof of theorem 4.5 will be completed in 5.2.

\demo{\bf{16.5. Ergodic cocycle and measurable transfer}}
Let $G$ be a locally compact group, $H$ a closed subgroup and
$X = G/H$.
Suppose $G$ has an invariant finite Borel measure
$\mu$
on $X$. Suppose we have a measurable section
$\delta : X \rightarrow G$.
For any
$g \in G$
and
$x \in X$
we define
$\lambda (g,x) \in H$
as unique element such that

$$s(gx) = gs(x) \lambda (g,x)$$
This defines a map of groups

$$G \overset {\lambda}\to \longrightarrow H^X$$
Now, $G$ acts on
$H^X$
by changing the argument. The map
$\lambda$
is well-known to be a {\it cocycle} for first non-abelian cohomology. Suppose
$f(h_1, \ldots, h_n)$
is a measurable $n$-cocycle on $H$. Then the composition
$f \circ \lambda : G^{(n)} \rightarrow \R$
is a measurable cocycle valued in the $G$-module of measurable functions on $X$. If $f$ is bounded, so is
$f \circ \lambda$.
In this case,
$\int_X f \cdot \lambda$
is a real bounded $n$-cocycle on $G$ and we have a well-defined map

$$H^n_b(H,\R) \overset t \to \longrightarrow H^n_b (G,\R)$$
Moreover, the composition
$H^n_b(G,\R) \overset {res} \to \longrightarrow H^n_b(H,\R) \overset t \to \longrightarrow H^n_b (G,\R)$
is a multiplication by
$\text{Vol} \, X$.
\enddemo

The proof of all this facts is easily adopted from [Gr1].

\demo{16.5.2}
As an immediate corollary, we state:
\enddemo

\proclaim{Theorem (16.5.2)}
Let
$\Gamma$
be a lattice in either 1)
$SO(n,1)$
or 2)
$SU(n,1)$
or 3)
$SU_{\H}(n,1)$.
Then in case 1)
$H^n_b(\Gamma) \ne 0$;
in case 2)
$H^{2\kappa}_b(\Gamma, \R) \ne 0$
for all
$1 \le \kappa \le n$;
in case 3)
$H^{4\kappa}_b(\Gamma,\R)$
for all
$1 \le \kappa \le n$.
\endproclaim

\demo{Proof}
Let us prove 2), since the rest is similar.
$SU(n,1)$
acts isometrically on the complex ball
$B^n$.
The K\"ahler form 
$\omega$
defined, as in 4.1, a class
$\omega$
in
$H^2_{b,cont}(SU(n,1),\R)$.
It is nontrivial since for any cocompact
$\Gamma$
the restriction on
$H^2(\Gamma)$
gives the K\"ahler class of
$B/\Gamma$.
Now for any
$\Gamma$,
the restriction on
$H^2_b(\Gamma,\R)$
must be nontrivial, other- wise
$\text{Vol} (SU(n,1)/\Gamma) \cdot \omega =0$
even as a class in
$H^2(SU^{\delta} (n,1),\R)$.
\enddemo

\demo{Proof of the theorem 16.4.5}
We need only to handle the case
$H_2(\Gamma,\R) =0$.
Then the restriction of the just defined class in
$H^2_b(\Gamma,\R)$
on
$H^2(\Gamma, \R)$
is zero, so by 18.4.2 we have a function $f$ with desired properties.
\enddemo

\bf{Appendix 17. Homology of ramified coverings} 

This appendix pursues a goal to develop further the three-manifold class field theory, started in [Re1]. In that paper we proved structure theorem for homology of a $p$-fold unramified covering $N$ of a 3-manifold $M$. In particular, we showed that
$H_1(N)_{(p)}$
as a 
$C_p$
-module is either a direct sum of closed modules, or a direct sum of exactly one open module and several closed modules. The terminology refers to Nazarova-Royter classification of finite modules over
$C_p$,
described in [Re1], 11.1. Moreover,
$H_1(N)_{(p)}$
is a sum of closed modules if and only if the covering is anisotropic [Re1], 11.4.

In the present paper we give a structure theorem for homology of a ramified covering. According to the ideology of arithmetic topology [KR] this corresponds to classical problems in number theory which have been studied since Gauss. The main result of the
 first part of this paper is the following:

\proclaim{Theorem A7}
Let
$N \rightarrow M$
be a covering of rational homology spheres ramified along a link with $s$ components. Then as a 
$C_p$-module,
$H_1(N)_{(p)}$
is a direct sum of exactly
$(s-1)$
open modules and some number of closed modules.
\endproclaim

If $M$ is a $p$-homology sphere, then we show that
$H_1(M)_{(p)} = \sum^{s-1}_{i=1} \Cal O/(1-\zeta)^{n_i}$,
where
$\Cal O=\Z[\zeta]/(1+\zeta + \ldots + \zeta^{p-1})$
is the ring of cyclotomic integers. The determination of numbers
$n_i$
is parallel to a classical problem in number theory [Red]. As a first step, we prove the following

\proclaim{Theorem A8}
Let
$p=2$
and
$s=2$,
so that the link is
$K_1 \cup K_2$
\roster
\item"(a)" if the linking number
$\text{link}\, (K_1,K_2)$
is odd, then
$H_1(M)_{(2)} = \Z/2\Z$.
\item"(b)" if the linking number
$\text{link}\, (K_1,K_2)$
is even, then either
$b_1(N) > 0$
or
$H_1(N)_{(2)} = \Z/2^{\ell}\Z$
with
$\ell \ge 2$.
\endroster
\endproclaim

We remark that some work in the direction of determination of
$n_i$
is being done by Kapranov [K] by entirely different methods.

Our methods in the paper develop those of [Re1]: spectral sequences for equivariant cohomology, transfer, reciprocity lemma, etc. Starting from secton 4 we make an intensive use of embedded nonoriented surfaces. Any 3-manifold with rich fundamental group 
(condition (R) of [Re1], 9.1) has a finite covering with embedded nonoriented surface, so this is a ``generic'' case.

\bf{17.1. Two spectral sequences of equivariant cohomology} 

\demo{17.1.1}
Consider a space $X$ with an action of a cyclic group
$C_p$.
Let
$E \rightarrow B$
be a classifying fibration for
$C_p$.
Recall that the equivariant cohomology
$H^*_e(X)$
is defined by
$H^*_e(X) = H^*(X \times E/C_p,\Z)$.
Since
$X \times E/C_p$
is a fibration over $B$ with fiber $X$, one arrives [B] to the first spectral sequence, converging to
$H^*_e(X)$
with
$E^2$
term
$H^i(C_p,H^j(X))$.
On the other hand, let
$Y = X/C_p$.
There is a sheaf
$\Cal F_i$
on $Y$ with a stalk
$(\Cal F_i)_y = H^i(G_y,\Z)$.
Here
$G_y$
is a stabilizer of any point
$x \in X$
over $y$. The second spectral sequence, converging to
$H^*_e(X)$
has
$E^2$-term
$H^i(X,\Cal F_j)$,
see [Q]. Observe that
$\Cal F_0$
is the constant sheaf
$\Z$.
Let
$F \subset X$
be the fixed point set for the 
$C_p$
-action; we may think of $F$ as a subspace of $Y$. Then
$\Cal F_i, i \ge 1$
all have support in $F$ and
$\Cal F_i|F$
is the constant sheaf
$\Z/p\Z$
for $i$ even and zero for $i$ odd.
\enddemo

\demo{17.1.2}
Now let $X$ be a closed oriented three-manifold. If
$p=2$
we assume that the action of
$C_p$
is orientation-preserving. Then the fixed point set $F$ consists of finite number
$\kappa$
of circles. It follows immediately that the second spectral sequence will have the
$E^2$
-term
\TagsOnLeft

$$\CD \ldots @.\qquad @. \ldots\\
0 @. \qquad @. 0 \\
\Z^{\kappa}_p @. \qquad @. \Z^{\kappa}_p\\
0 @. \qquad @. 0 \\
\Z^{\kappa}_p @. \qquad @. \Z^{\kappa}_p \\
0 @. \qquad @. 0 @. \qquad @. 0 @. \qquad @. 0\\
\Z @. \qquad @. H^1(Y,\Z) @. \qquad @. H^2(Y,\Z) @. \qquad @. \Z \endCD \tag II$$
where we silently used the fact that
$Y = X/C_p$
is a closed manifold. The first spectral sequence will have the following
$E^2$
-term:

$$\CD \Z @. \qquad @. 0 @. \qquad @. \Z_p @. \qquad @. 0 @. \qquad \Z_p @. \qquad @. \ldots \\
H^i(C_p, H^2(X)) \\
H^i(C_p,H^1(X)) \\
\Z @. \qquad @. 0 @. \qquad @. \Z_p @. \qquad @. 0 @.\qquad \Z_p @.\qquad @. 0 @.\qquad @. \ldots 
\endCD \tag I$$
>From now on we assume that $X$ is a rational homology sphere, that is,
$H_1(X,\Z)$
is finite. Then it follows from I that
$\dim_{\F_p} H^2_e(X) = 1 + \dim_{\F_p}H^0(C_p,W)$.
Here
$W = H_1(X,\Z)_{(p)}$
and for a finite abelian $p$-group
$A, \dim_{\F_p} A = \log_p |A|$.
Moreover,
$H^2(X,\Z)_{(p)} \approx \hat W \approx W$,
see [Re1], 5.1.

Now, one deduces from II that
$\dim H^2_e(X) = \kappa + \dim H_1(Y)_{(p)}$.
As a result, the following proposition follows.
\enddemo

\proclaim{Proposition (17.1.2)}
$\dim H^0(C_p,W) = \dim H_1(Y)_{(p)} + \kappa -1$.
\endproclaim

\bf{17.2. The transfer map for ramified coverings and structure theorems} 

>From now on we relabel $X$ by $N$ and $Y$ by $M$, so
$N \overset {\pi} \to \longrightarrow M$
is a normal covering, ramified over a link with $k$ components. Let
$K_i, \, 1\le i \le \kappa$
be these components,
$K_i \subset M$.
We wish to define a map

$$t : H_1 (M,\Z) \rightarrow H_1(N,\Z)$$
is the following fashion. For $L$ a link disjoint from
$\cup K_i$
consider
$\pi^{-1}(L) \subset M$.
If
$L_1, L_2$
are two homologous links in $M$, there is a surface
$F \subset M$
with
$\partial F = L_1-L_2$.
We may assume that $F$ is transversal to
$\cup K_i$.
Moreover, we can choose a tubular neighbourhood $U$ of
$\cup K_i$
with coordinates
$(\tau, z)$
such that
$\pi|\pi^{-1}(U)$
takes the form
$(\tau, w) \mapsto (\tau, w^p)$.
Here $w$ runs in a small disc in
$\C$.
We can further arrange that
$F \cap U$
is a union of z-discs, so that
$\pi^{-1}(F)$
is a smooth surface and
$\pi : \pi^{-1}(F) \rightarrow F$
is a ramified covering. It is clear now that
$\pi^{-1}(L_1)$
is homologous to
$\pi^{-1}(L_2)$,
so $t$ is well-defined.

\proclaim{Proposition (19.2.1)}
$\pi t = p \cdot id, t\pi = 1+\zeta + \ldots + \zeta^{p-1}$,
where
$\zeta : H_1(N) \rightarrow H_1(N)$
is a generator of the
$C_p$
-action.
\endproclaim

The proof is obvious.

If
$M,N$
are rational homology spheres, then we have the linking form in
$H_1$
(see [Re1], 5.1). The following important statement is parallel to [Re1], 5.2.

\proclaim{17.2.2 Proposition (Reciprocity)}
$t = \pi^*$.
In words, for
$z \in H_1(M)$
and
$v \in H_1(N)$,

$$(\pi v,z)_M = (v,tz)_N$$
\endproclaim

\demo{Proof}
Let $N$ be such that
$N \cdot [z] =0$
in
$H_1(M)$.
Let
$\gamma, \delta$
be knots representing
$z,v$. Let $F$ be a 2-cycle in $M$ with boundary
$N \cdot \gamma$.
Acting as above, we can achieve that
$\gamma$
is disjoint from the ramification locus and
$\pi^{-1}(F)$
is a two-cycle with boundary
$N(\pi^{-1}(\gamma))$.
Then the proof goes as in [Re1].
\enddemo

\proclaim{17.2.3 Theorem (Injectivity of transfer)}
Suppose
$\kappa \ne 0$,
that is,
$\pi$
is indeed ramified. Then
$t : H_1(M)_{(p)} \rightarrow H_1(N)_{(p)}$
is injective.
\endproclaim

Observe that this is never the case if
$\pi$
is unramified ([Re1]).

\demo{Proof}
It is enough to show that if
$z \in H_1(M)_{(p)}$
has period $p$, then
$tz \ne 0$. Consider the cohomology class
$\lambda = (\cdot, z) : H_1(M) \rightarrow \F_p$.
It defines an unramified covering
$Q \rightarrow M$.
If
$\lambda$
vanished on the image of
$H_1(N)$
in
$H_1(M)$,
we would have a lift

$$\CD N @>>> Q \\
@AAA \swarrow\\
M \endCD$$
which is clearly impossible. So
$(\pi v,z) \ne 0$
for some $v$,
hence
$(tz,v) \ne 0$
by 2.2.
\enddemo

\demo{17.2.4 Corollary}
$\pi : H_1(N)_{(p)} \rightarrow H_1(M)_{(p)}$
is surjective.
\enddemo

We now arrive to the following very important result.

\proclaim{19.2.5 Theorem}
Suppose
$\kappa \ne 0$.
Then for all
$i \ge 1$,

$$\dim_{\F_p} H^i(C_p,W) = \kappa -1$$
Recall that if
$\kappa =0$
then
$\dim H^i(C_p,W)$
is either 0 or 1 by [Re1], 2.1.
\endproclaim

\demo{Proof}
Since
$\dim H^1(C_p,W) = \dim H^2(C_p,W)$
[$\quad$], it is enough to show this for
$k=2$.
By 1.2,
$\dim H^0 (C_p,W) = \dim H_1(M)_{(p)} + \kappa -1$.
Moreover,
$\text{Im} (1+\zeta + \ldots + \zeta^{p-1}) = \text{Im} \, t$
by 2.1 and 2.4. On the other hand,
$\text{Im} \, t \approx H_1 (M)_{(p)}$
by 2.3, so
$\dim \text{Im} (1+\zeta + \ldots + \zeta^{p-1}) = \dim H_1(M)_{(p)}$.
It follows that
$\dim H^2(C_p,W) = \dim H^0(C_p,W)/\text{Im} (1+\zeta + \ldots + \zeta^{p-1}) = \kappa -1$.
\hfill Q.E.D.
\enddemo

\demo{17.2.6 Corollary}
Suppose we have two isotopic actions of
$C_p$
on a rational homology sphere $N$ with fixed point sets consisting of
$\kappa$
and
$\ell$
circles respectively. If
$\kappa \ge 2$,
then
$\kappa = \ell$.
\enddemo

For the following corollary the reader needs to know the Nazarova-Royter theory explained in [Re1].

\proclaim{A.7 Theorem (Structure theorem)}
Suppose
$k \ne 0$.
Then as a
$C_p$
-module,
$W = H_1(M)_{(p)}$
is a direct sum exactly
$(\kappa -1)$
open modules and some number of closed modules.
\endproclaim

\proclaim{Corollary (19.2.8)}
Suppose $M$ is a $p$-homology sphere. Then
$W = \bigoplus^{n-1}_{i=1} \Cal O/(1-\zeta)^{\ell_i}$.
Here
$\Cal O$
is the ring of integers in the cyclotomic field.
\endproclaim

\demo{Proof}
Since
$H_1(M)_{(p)} =0$
we deduce that
$1 + \zeta + \ldots + \zeta^{p-1} = t \cdot \pi =0$,
so $W$ is an
$\Cal O$
-module, hence a sum of
$\Cal O/(1-\zeta)^{\ell_i}$.
There should be
$(\kappa -1)$
-summands by A.7.
\enddemo

For
$p=2$ 
we specialize 19.2.8 to:

\proclaim{Corollary (17.2.9)}
Suppose
$p=2$
and $M$ is a 2-homology sphere. Then
$W = \bigoplus^{\kappa-1}_{i=1} \Z/2^{\ell_i} \Z$
and
$C_2$
acts by multiplication by 
$(-1)$.
\endproclaim

\bf{17.3. Ramification over 2-component links, I}

\demo{17.3.1} We start with a following general observation.
\enddemo

\proclaim{Proposition (17.3.1)}
Suppose
$N \rightarrow M$
is as above. Then:
\roster
\item"(a)" $\Sigma [K_i] \in p \cdot H_1(M,\Z)$
\item"(b)" if
$\Sigma [K_i] \in np \cdot H_1(M,\Z)$
then in
$H_1(N,\Z), \Sigma[K_i] \in n \cdot H_1(N,\Z)$
\endroster
\endproclaim

\demo{Proof} 
(a) The $m$-ramified covering 
$N \rightarrow M$
ramified over
$\cup K_i$
exists if and only if
$\Sigma [K_i] \in m \cdot H_1(M;\Z)$
(see [H]). This induces (a). If
$\Sigma [K_i] \in np \cdot H_1(M,\Z)$,
then there is a $np$
-ramified covering, of $M$, ramified over
$\cup K_i$,
therefore there exists a $n$-ramified covering over $N$, so
$\Sigma [K_i[ \in n \cdot H_1(N,\Z)$.
\hfill Q.E.D.
\enddemo

\demo{Corollary (17.3.2)}
If
$\Sigma[K_i] =0$
in
$H_1(M,\Z)$,
then the same holds in
$H_1(N,\Z)$.
\enddemo

\demo{17.3.2}
Now let $M$ be a 2-homology sphere, and
$K_1 \cup K_2$
be a two-component link in $M$. Let
$N \rightarrow M$
be a 2-fold ramified covering, ramified over
$K_1 \cup K_2$.
We know by 2.8, that either
$b_1(N) > 0$,
or
$H_1(N)_{(2)} = \Z/2^{\ell} \Z$.
\enddemo

\demo{Example (17.3.2.1)}
Let
$M = S^3$
and
$K_1, K_2$
be two linked great circles. Then
$N \approx \R P^3$,
in particular
$H_1(N)_{(2)} = \Z_2$.
\enddemo

\demo{Example (17.3.2.2)}
Let
$M = S^3$
and
$K_1, K_2$
be two separated knots (two knots lying in disjoint balls). Then
$b_1(N) = 1$.
\enddemo

These examples are the motivation for the following structure theorem.

\proclaim{Theorem (17.3.2)}
Let $d$ be a linking number 
$\text{mod} \, 2$
of
$(K_1, K_2)$.
Then
\roster
\item"(a)" if
$d=1$
then
$b_1(N) =0$
and
$H_1(N)_{(2)} = \Z_2$.
\item"(b)" if
$d=0$
then either
$b_1(N) > 0$
or
$b_1(N) =0$
and
$H_1(N)_{(2)} = \Z/2^{\ell}\Z, \ell \ge 2$.
\endroster
\endproclaim

The proof of this theorem will be completed in section 5.

\demo{Corollary (17.3.3)}
Let $N$ be a three-manifold such that
$b_1(N) =0$
and
$H_1(N)_{(2)} = \Z_2$.
Let
$\zeta$
be an involution of $N$. Then either
$\zeta$
is free, or
$M = N/\zeta$
is a 2-homology sphere,
$\text{Fix} (\zeta) = K_1 \cup K_2$
and the linking number of
$(K_1,K_2)$
in $M$ is odd.
\enddemo

\demo{Proof}
Suppose
$\zeta$
is not free. Since
$H^1(C_2,H_1(N)_{(2)}) = \Z_2$,
we know by 2.7 that
$\text{Fix}(\zeta)$
consists of two components, say
$K_1$
and
$K_2$.
By 1.2, 
$\dim H_1(M)_{(2)} =0$
so $M$ is a 2-homology sphere. Finally, the linking number of
$(K_1,K_2)$
may not be even by 3.2.
\enddemo

\bf{17.4. Proof of the Theorem 17.3.2} 

Let
$F_i$
be a (possibly nonorientable) surface with boundary
$K_i$.
We can always arrange that
$F_1$
interesect a tubular neighbourhood of
$K_2$
by
z-discs, and the same for
$F_2$.
Then
$\pi^{-1}(F_i)$
is a smooth surface {\it without boundary}. Moreover
$\pi : \pi^{-1}(F_1) - K_1 \rightarrow F_1 - K_1$
is a double covering ramified over
$F_1 \cap K_2$.
It follows that
$\chi (\pi^{-1}(F_i)) = \text{link} \, (K_1,K_2) (\text{mod} \, 2)$.
Assume
$\text{link} \, (K_1,K_2) =1 (\text{mod} \, 2)$,
then
$\chi(\pi^{-1}(F_i)) =1 (\text{mod} \, 2)$
in particular,
$\pi^{-1}(F_i)$
is nonorientable. Moreover,
$[K_1] \cap [\pi^{-1}(F_2)] =1$.
We know from 2.8 that
$\zeta =-1$
on
$H_1(N,\Z)$.
Suppose
$b_1(N) > 0$.
Then
$H^0(C_2,H^1(N,\Z)) =0$,
so that the differential
$d_2 : H^0(C_2,H^1(N,\Z)) \rightarrow \Z_2$
in the  spectral sequence I of 1.2 is zero. Therefore
$\dim H^2_e(N) = 1 +\dim H^0(C_2,H^2(N)) + \dim H^1(C_2,H^1(N))$.
>From the spectral sequence II we find
$\dim H^2_e(N) =2$.
So
$\dim H^0(C_2,H^2(N)) + \dim H^1(C_2,H^1(N)) =1$.
However, if
$H_1(N) = \Z^m \oplus \Z/2^{\ell_1} \Z \oplus \ldots \oplus \Z/2^{\ell_s} \Z \oplus \Sigma \Z/p^{r_i}_i \Z, p_i$
odd, then
$\dim H^0(C_2,H^2(N)) = s$
and
$\dim H^1(C_2,H^1(N)) = m$.
Hence
$m+s=1$,
so either
$H_1(N) = \Z \oplus \Sigma \Z/p^{r_i}_i \Z$,
or
$H_1(N) = \Z/2^{\ell} \Z \oplus \Sigma \Z/p^{r_i}_i \Z$.
In the first case the natural map
$H^1(N,\Z) \rightarrow H^1(N,\F_2)$
is surjective, therefore the multiplication in
$H^1(N,\F_2)$
is zero. However, since
$\chi(\pi^{-1}(F_i))$
is odd,
$\pi^{-1}(F_i) \cap \pi^{-1} (F_i) \cap \pi^{-1}(F_i) \ne 0$.
So the first case is impossibe, and
$b_1(N) =0$
and
$H_1(N)_{(2)} =\Z/2^{\ell}\Z$.
By [Re1], 1.6.1,
$\ell =1$,
which proves (a).

Now let
$\text{link}(K_1,K_2) =0 (\text{mod} \, 2)$
and
$b_1(N) =0$.
Then
$H_1(N)_{(2)} = \Z/\Z^{\ell} \Z$
by 2.8 and we need to show that
$\ell \ne 1$.
Suppose the opposite, that is,
$H_1(N)_{(2)} = \Z/2\Z$.
If
$\pi^{-1}(F_1)$
were nonorientable, it would represent a generator in
$H_2(N,\F_2)$
by 4.1. Since
$\chi(\pi^{-1}(F_1))$
is even, we would have
$\lambda^3 =0$
for the dual class. However,
$\lambda^2= \beta(\lambda) \ne 0$,
hence
$\lambda^3 \ne 0$, in contradiction. So
$\pi^{-1}(F_1)$
is orientable and separates $N$.

By assumption,
$K_1$
is homologically zero. That means we can choose
$F_1$
to be orientable. Now we claim that
$K_1$
as a curve in
$\pi^{-1}(F_1)$,
separates
$\pi^{-1}(F_1)$.
Indeed, the involution
$\zeta : N \rightarrow N$
restricted on
$\pi^{-1}(F_1)$
obviously reserves the orientation of
$\pi^{-1}(F_1)$.
The same holds for
$\pi^{-1}(F_1)- K_1$.
On the other hand,
$\zeta$
preserves the orientation of
$\pi^{-1}(F_1)- K_1$,
induced from
$F_1- K_1$.
So
$\pi^{-1}(F_1)- K_1$
is disconnected, and
$\zeta$
permutes the components. It follows that
$K_2 \cap F_1 = \emptyset$.
But then
$K_2$
should lie in one of the components of
$N - \pi^{-1}(F_1)$,
which are permuted by
$\zeta$,
which is impossible. So
$\ell \ne 1$.

\bf{17.5. The butterfly diagram for Hilden-Montesinos knots} 

In this section we will describe a very interesting class of examples, illustrating the situation in Theorem 3.2 (b). Let $M$ be a rational homology sphere with
$H_1(M)_{(2)} = \Z/2\Z$.
By [H], [M] there is a knot
$K \subset S^3$
such that $M$ is an irregular 3-covering of
$S^3$
ramified over $K$. It follows that there is a double ramified covering $F$ of
$S^3$
ramified over $K$, and a double ramified covering $N$ over $M$, ramified over $K$, such that $N$ is a Galois unramified
$C_3$
-covering of $F$:

\midspace{4cm}
$   $\newline
Notice that $F$ is a 2-homology sphere [BZ]. Assume
$b_1(N) =0$.
Then
$H^i(C_2,H_2(N)_{(2)}) =0, (i \ge 1)$
by 2.5. Moreover, 
$H^0(C_2,H_2(N)) = \Z_2$
by 1.2. It follows from [Re1], 11, that either
$H_1(N)_{(2)} = \Z/2\Z \oplus \Z/2 \Z$,
or
$H_1(N)_{(2)} = \Z/2^{\ell} \Z, \ell \ge 3$
and
$C_2$
acts by multiplication by
$\pm 1 + 2^{\ell-1}$.
However,
$C_3$
acts on
$H_1(N)_{(2)}$
fixed point free, therefore
$H_1(N)_{(2)} = \Z/2\Z \oplus \Z/2\Z$
and the linking form is hyperbolic [Re1], 12.5. It follows from [Re1], 12.5, that either $N$ has a quaternionic covering with positive
$b_1$,
or
$N = R/Q_8$
where $R$ is a 2-homology sphere. So either $F$ has a Galois covering of order
$\le 24$
with positive
$b_1$, or we have the following diagram:

\midspace{6cm}
$   $\newline
Now, let
$Q \rightarrow M$
be a unique unramified 
$C_2$
-covering. Then we can realize $S$ as a pullback:

\midspace{4 cm}
$   $\newline
in particular,
$b_1(Q) =0$
so
$H_1(Q)_{(2)} =0$
by [Re1], 7.2. Since
$H_1(S)_{(2)} = \Z/2\Z$
it follows that
$S \rightarrow Q$
is ramified over a 2-component link
$K_1 \cup K_2$
by 2.8. Moreover, link
$(K_1,K_2)$
is even by Theorem 3.2.

Let
$\zeta, \eta$
be two commuting involutions of $S$, corresponding to coverings
$S \rightarrow N$
and
$S \rightarrow Q$,
so that
$\zeta$
is free and
$\eta$ 
is not. By [Re1], 12.5, the action of
$\zeta$
on
$H_1(S)_{(2)} = \Z/4\Z$
is multiplication by 
$(-1)$.
By 2.8, 
$\eta$
acts also by multiplication by
$(-1)$.
The free involution
$\zeta\eta$
will act identically so by [Re1], 11, the manifold
$P = S/\zeta \eta$
has homology
$H_1(P)_{(2)} = \Z/8 \Z$.
Summing up, we have a diagram

\newpage
$   $
or else $F$ is virtually Haken.

\bf{17.6. An ierarchy of 2-homology spheres} 

In this chapter we will define, for a link
$(K_1,K_2)$
in a 2-homology sphere
$M_0$,
a very interesting sequence of three-manifolds, associated to the link. The definition is based on the theory of developed in chapter 3-6.

So let
$(K_1,K_2)$
be a link in $M$ satisfying the following assumptions:
\roster
\item"(i)" $K_i$
is homologically trivial
\item"(ii)" the (usual) linking number
$\text{link}\,(K_1,K_2)$
is odd.
\endroster

\proclaim{Theorem (17.6.1)}
Let $N$ be a covering of
$M_0$
ramified over
$K_1 \cup K_2$.
Let
$M_1$
be the unique double unramified covering of $N$ and let
$K^{(1)}_1$
and
$K^{(1)}_2$
be preimages of
$K_1,K_2$
in
$M_2$.
Then either
$b_1(M_1) > 0$
or
\roster
\item"(a)" $M_2$
is a 2-homology sphere
\item"(b)" $K^{(1)}_i$
are connected
\item"(c)" $K^{(1)}_i$
are homologically trivial
\item"(d)" $\text{link}(K^{(1)}_1, K^{(1)}_2) = \text{link}(K_1,K_2)$
\endroster
\endproclaim

\demo{Proof} Suppose
$b_1(M_1) =0$,
then
$M_1$
is a 2-homology sphere by 3.2 and [Re1], 7. Let
$F_i$
be an oriented Seifert surface of
$K_i$,
intersecting the other component by z-discs. If
$\pi : N \rightarrow M$
is the covering, we call
$Q_i = \pi^{-1}(F_i)$.
We know that
$[Q]$
is the generator in
$H_2(N,\F_2)$
and
$[Q_1] \cap [K_2]=1$.
It follows that
$K^{(1)}_i$
are connected. Let
$\delta : M_1 \rightarrow N$
be the covering and let
$S_i = \delta^{-1}(Q_i)$.
By 4.1, if
$\lambda$
is the dual class to
$[Q_i]$
in
$H^1(N,\F_2)$,
then
$\lambda|Q_i = w_1(Q_i)$,
so that
$S_i$
are oriented. Let
$\zeta$
be the Galois involution of $N$, then
$\text{Fix} (\zeta ,Q_i) = K_i$.
Now, 
$\zeta$
lifts to an involution of
$S_i$
such that
$\text{Fix}(\zeta, S_i) = K^{(1)}_i$.
Moreover, if
$\eta$
is the free Galois involution of
$M_1$,
then
$\zeta$
and
$\eta$
commute. So we have a 
$C_2 \times C_2$
action on
$S_i$,
which is free on
$S_i - K_i$.
Since
$S_i - K_i /C_2 \times C_2 \approx F_i - K_i$,
this action preserves an orientation of
$S_i - K_i$.
On the other hand,
$\eta$
reverses the orientation of
$S_i$.
It follows that
$S_i - K_i$
is disconnected, so
$K_i$
bounds an oriented surface, which proves (c). Now (d) follows from above by direct computation.

In this way we define an ierarchy

$$\ldots \rightarrow N_1 \rightarrow M_1 \rightarrow N_0 \rightarrow M_0$$
where
$M_i$
are 2-homology spheres with a link
$(K^{(i)}_1, K^{(i)}_2)$,
and
$H_1(N_i)_{(2)} = \Z/2Z$.
\enddemo

Permanent address: Department of Mathematical Sciences, University of
Durham, Durham, UK.

Current address (till
June, 1998)

CNRS and Institut Fourier, BP 74, 38402 Saint Martin d Heres, France.

email: reznikov\@ mozart.ujf-grenoble.fr

%\input biblio...
%Musterdokument fuer Literaturangaben: biblio.tex

\end